\documentclass[12pt,a4paper]{amsart}

\usepackage{tikz-cd}
\usepackage[margin=1in]{geometry}  % set the margins to 1in on all sides
\usepackage{graphicx,hyperref}              % to include figures
\usepackage{amsmath,euscript}               % great math stuff
\usepackage{amsfonts}              % for blackboard bold, etc
\usepackage{amssymb,,amscd,epsf,amsthm, epsfig}                % better theorem environments
\usepackage[toc,page]{appendix}

\newtheorem{thm}{Theorem}[section]

\newtheorem{prop}[thm]{Proposition}
\newtheorem{cor}[thm]{Corollary}
\newtheorem{lemma}[thm]{Lemma}

\newtheorem{obs}[thm]{Observation}
\theoremstyle{definition}
\newtheorem{rem}[thm]{Remark}
\newtheorem{defn}[thm]{Definition}

\input xy
\xyoption{all}

\newcommand{\Z}{\mathbb{Z}}
\newcommand{\LL}{\mathbb{L}}

\newcommand{\calM}{\mathcal{M}}

\newcommand{\calW}{\mathcal{W}}
\newcommand{\calR}{\mathcal{R}}

\newcommand{\MM}{\mathbb{M}}

\newcommand{\DD}{\EuScript D}

\newcommand{\pf}{{\it Proof.}\hspace{2ex}}

\newcommand{\epf}{\hspace*{\fill}\mbox{$\halmos$}}
\newcommand{\halmos}{\rule{1ex}{1.4ex}}

\newcommand{\Hom}{\mathop{\mathrm{Hom}}\nolimits}

\newcommand{\Imm}{\mathop{\mathrm{Im}}\nolimits}

\newcommand{\rad}{\mathop{\mathrm{rad}}\nolimits}

\newcommand{\id}{\mathop{\mathrm{id}}\nolimits}

\newcommand{\op}{\mathop{\mathrm{op}}\nolimits}

\newcommand{\Ker}{\mathop{\mathrm{Ker}}\nolimits}
\newcommand{\Tor}{\mathop{\mathrm{Tor}}\nolimits}

\newcommand{\sg}{\mathop{\mathrm{sg}}\nolimits}

\newcommand{\Sl}{\mathop{\mathfrak{sl}}\nolimits}

\newcommand{\HH}{\mathop{\mathrm{HH}}\nolimits}

\begin{document}

%\nocite{*}

\title{Tate-Hochschild Cohomology of Radical Square Zero Algebras}
\date{}

\thanks{}

\date{\today}

\author{Zhengfang Wang}
\address{Zhengfang Wang\\ Beijing International Center for Mathematical Research, Peking University, No. 5 yiheyuan Road Haidian District, 100871 Beijing, China}
\address{Universit\'e Paris Diderot-Paris 7, Institut
de Math\'ematiques de Jussieu-Paris Rive Gauche CNRS UMR 7586, B\^atiment Sophie Germain, Case 7012,
75205 Paris Cedex 13, France}

\email{wangzhengfang@bicmr.pku.edu.cn}

\maketitle
\begin{abstract}
For algebras with radical square zero, we give a combinatorial description to the Tate-Hochschild cohomology. We compute the Gerstenhaber algebra structure on the Tate-Hochschild cohomology for some classes of such algebras.
\end{abstract}

\section{Introduction}
Following Buchweitz \cite{Bu} and Orlov \cite{Orl} the singularity category $\DD_{\sg}(A)$ of a Noetherian  algebra $A$ is the localization of the bounded derived category of $A$ at the full subcategory of compact objects. For more details see e.g. [Zim].  In a recent paper \cite{Wang}, we studied the Tate-Hochschild cohomology ring $\HH_{\sg}^*(A, A)$ of a Noetherian  algebra $A$ over a commutative ring $k$ such that $A$ is projective as a $k$-module.   By definition,  $$\HH_{\sg}^i(A, A):=\Hom_{\DD_{\sg}(A\otimes_k A^{\op})}(A, A[i]),$$ where we denote as usual by $[i]$ the $i$-th iterate of the suspension functor for $i\in\Z$. In loc. cit.  we also proved that the Tate-Hochschild cohomology ring $\HH^*_{\sg}(A, A)$
is a Gerstenhaber algebra whose cup product is given by the Yoneda product in $\DD_{\sg}(A\otimes A^{\op})$ (cf. \cite[Proposition 4.7]{Wang}). %The Tate-Hochschild cohomology group $\HH_{\sg}^i(A, A)$ can be computed by the bar resolution of $A$, roughly speaking, it is realized as the colimit of an inductive system associated to the bar resolution (cf. \cite[Proposition 3.1]{Wang}). 
By the very recent work of Keller \cite{Kel}, the Tate-Hochschild cohomology of an algebra $A$ is isomorphic, as graded algebras,  to the Hochschild cohomology of the dg singularity category (i.e. the canonical dg enhancement of the singularity category) of $A$. 
For further studies on the Tate-Hochschild cohomology, we refer to \cite{BeJo, EuSc, Ngu, RiWa, Wang}.

A radical square zero algebra is a finite dimensional algebra over a field $k$ such that
the square of its Jacobson radical is zero.
 From Gabriel's theorem it follows that every radical square zero algebra over an algebraically closed field $k$ is Morita equivalent to a path algebra modulo relations of the form
$kQ/\langle Q_2\rangle$, where $Q$ is a finite quiver and $Q_2$ is the set of all paths in $Q$ of length 2. In this paper, by a {\it radical square zero algebra} we mean an algebra $kQ/\langle Q_2\rangle$ for a finite quiver $Q$ over a field $k$ (not necessarily algebraically closed). 

Following the approach developed by Cibils and S\'{a}nchez (cf. \cite{Cib1, Cib2, Cib3, San, San2}), we study the Tate-Hochschild cohomology of radical square zero algebras. In particular, we will compute  the Gerstenhaber algebra structure on the Tate-Hochschild cohomology for some classes of such algebras. Recall that in \cite{Cib1}, Cibils gave a very combinatorial description for the Hochschild cohomology ring $\HH^*(A, A)$ in the case of a radical square zero algebra $A$. Following this construction, in \cite{San, San2} S\'{a}nchez studied the Lie module structure on the Hochschild cohomology groups $\HH^*(A, A)$
over the Lie algebra $\HH^1(A, A)$. Inspired by these work, we give a combinatorial descriptions to the Tate-Hochschild cohmology $\HH_{\sg}^*(A, A)$ of radical square zero algebras $A$. 
%We also provide examples of algebras  whose Hochschild cohomology (or Tate-Hochschild cohomology)
%do not admit a BV algebra structure (cf. Remark \ref{rem-conter1}).

This paper is organized as follows. Section 2 is devoted to recalling some  notions and results in \cite{Cib1, Cib2, Cib3}. In Section 3, we prove the main  result (cf. Theorem \ref{prop-singular}) using the notion of  parallel paths appearing in \cite{Cib1}. This gives a very combinatorial description to the Tate-Hochschild cohomology of radical square zero algebras.  In Section 4, we give examples of the radical square zero algebras associated to $c$-crown quivers. The Gerstenhaber algebra structure on $\HH^*_{\sg}(A, A)$ will be computed explicitly (cf. Proposition \ref{prop-de-ger} and Theorem \ref{equ-bv00}). Section 5  deals with  radical square zero algebras $A$ associated with $r$-multiple loops quivers. We give a complete description on the Gerstenhaber algebra structure on $\HH_{\sg}^*(A, A)$ (cf. Corollary \ref{cor-last}). %This can be used to study the Hochschild cohomology of Leavitt path algebras.
%In Section 6, we give a prop interpretation for the Gerstenhaber bracket in the case of radical square zero algebras. This generalizes the construction in Section 5. In Appendix A, we use prop theory to construct a new Lie algebra structure on $\Gl_{\infty}(k)$.

%For simplicity, now let us assume that the base ring $k$  is a field of characteristic zero.
%In this paper, we  frequently use some notions on quivers without recalling. We refer to
%\cite{Cib1, Cib3, Boe, Mich, Zim} for details.

Motivated by results in the present paper and in \cite{ChYa}, we establish,  in an ongoing work \cite{ChLiWa}, an isomorphism of Gerstenhaber algebras from the Tate-Hochschild cohomology of $kQ/\langle Q_2\rangle)$ to the Hochschild cohomology of the Leavitt path algebra $L_k(Q)$ (cf. e.g. \cite{AbAr}) for a finite quiver $Q$. From this point of view, the Tate-Hochschild cohomology may control the deformation theory of singularity categories. This topic will be explored in a future research.     We refer to \cite{Che, ChYa, Smi} for the studies of singularity categories of radical square zero algebras. 

\section*{Acknowledgement} This work is a part of author's PhD thesis. He would like to thank his supervisor Alexander Zimmermann for introducing this interesting topic and for his many valuable suggestions for improvement. He also would like to thank Huafeng Zhang for many useful discussions during this project.  The author is indebted to Ragnar-Olaf Buchweitz for the constant support and encouragement. 

The author is very grateful to the referee for valuable suggestions and comments.  The author was partially supported by NSFC (No.11871071).

\section{Background and Notation}\label{section2}
In this section, we assume that $k$ is a field and $A$ is a finite dimensional $k$-algebra.  Recall that the Jacobson radical $\rad(A)$ is defined as the intersection of all maximal left ideals of $A$. Wedderburn-Malcev theorem says that there exists a semi-simple subalgebra $E$ of $A$ such that $A\cong E\ltimes \rad(A)$ as $k$-algebras.  

Put $r:=\rad(A)$.  By Lemma 2.1 in \cite{Cib1}, we have a projective resolution $\calR(A)$ of the $A$-$A$-bimodule $A$: 
  \begin{equation}\label{equ-proj-res}
    \xymatrix{
  \cdots\ar[r]^-{d_{i+1}} & A\otimes_E r^{\otimes_E i}\otimes_E A \ar[r]^-{d_i} & \cdots\ar[r]^-{d_2} & A\otimes_E r\otimes_E A\ar[r]^-{d_1} &
    A\otimes_E A \ar[r]^-{\epsilon} & A\ar[r] & 0
    }
  \end{equation}
  where $\epsilon(a\otimes b)=ab$
  and
  \begin{equation*}
    \begin{split}
      d_i(a\otimes x_{1, i}\otimes b)=&ax_1\otimes x_{2, i}\otimes b+\sum_{j=1}^i(-1)^ja\otimes x_{1, j-1}\otimes x_jx_{j+1}\otimes x_{j+2, i}\otimes b+\\
     & (-1)^ia\otimes x_{1, i-1}\otimes x_ib.
    \end{split}
  \end{equation*}
 Here for simplicity, we write $x_i\otimes \cdots \otimes x_j\in r^{\otimes j-i}$ as $x_{i, j}$ and the tensor product  $\otimes_E$ over $E$ as $\otimes$. 
 
  Let $M$ be an $A$-$A$-bimodule.  Then the Hochschild cohomology groups $\HH^i(A, M)$ are cohomology groups of the {\it Hochschild cochain complex} $C^*(r, M)$ (cf. \cite[Proposition 2.2]{Cib1}):
  \begin{equation}\label{equ-proj-res1}
    \xymatrix{
    0\ar[r] & M^E\ar[r]^-{\delta} & \Hom_{E-E}(r, M)\ar[r]^-{\delta}  & 
  \cdots\ar[r]^-{\delta} &  \Hom_{E-E}(r^{\otimes_E i}, M) \ar[r]^-{\delta} & \cdots
    }
  \end{equation}
  where  $M^E=\{m\in M\ |\ sm=ms  \mbox{ for all $s\in E$}\}.$ For $m\in M^E$ and $x\in r$, $$\delta(m)(x)=mx-xm$$  and for $f\in\Hom_{E-E}(r^{\otimes_Ei}, M)$,
  \begin{equation*}
    \begin{split}
      \delta(f)(x_{1, i+1})=x_1f(x_{2, i+1})+\sum_j^{i} (-1)^{j} f(x_{1, j-1} \otimes x_jx_{j+1}\otimes x_{j+2, i+1})
      +(-1)^{i+1} f(x_{1, i}) x_{i+1}.
    \end{split}
  \end{equation*}
\begin{rem}\label{rem-twoterms}
 Observe that if $A$ is a radical square zero algebra, then the differential $\delta(f)$ has only two terms,  namely, $\delta(f)(x_{1, i+1})=x_1f(x_{2, i+1})+(-1)^{i+1} f(x_{1, i})x_{i+1}.$ 
 \end{rem}
Let $Q=(Q_0, Q_1, s, t)$ be a finite (i.e. $|Q_0\cup Q_1|<+\infty$) and  connected quiver. For $n\in\Z_{\geq 0}$, we define $Q_n$ to be the set of all paths in $Q$ of length $n$. The trivial path is denoted by $e_i$ for a vertex $i$ in $Q_0$. We consider the path algebra $kQ$ of $Q$ over the field $k$. As a vector space, $kQ=\bigoplus_{n\in\Z_{\geq 0}} kQ_n.$ If  $p$ and $q$ are two paths in $Q$, then their product $pq$ is the concatenation of the paths $p$ and $q$ if $t(q)=s(p)$ and $0$ otherwise. The radical square zero algebra of $Q$ is defined as $A_Q:=kQ/\langle Q_2 \rangle,$ where $\langle Q_2\rangle$ is the two-sided ideal in $kQ$ generated by $Q_2,$ the set of paths of length 2. Note that $A_Q$ is a finite-dimensional $k$-algebra with the Jacobson radical $\rad(A_Q)=kQ_1$ and the Wedderbrun-Malcev decomposition is $A_Q=kQ_0\ltimes kQ_1$, where $E=kQ_0$ is a commutative semi-simple algebra isomorphic to $\bigoplus_{e_n\in Q_0} ke_n.$

%In general, for a finite-dimensional $k$-algebra $A$ (not necessarily, radical square zero algebra) with a Wedderburn-Malcev decomposition $A=E\oplus \rad(A)$, we have a projective resolution of the $A$-$A$-bimodule $A$. We abbreviate in the sequel$$a_{i, j}:=a_i\otimes a_{i+1}\otimes \cdots\otimes a_j$$ for $i<j$.

Next we will apply the projective resolution $\calR(A)$ in (\ref{equ-proj-res}) to compute the Tate-Hochschild cohomology rings of
radical square zero algebras. From now on,  we fix a finite and connected quiver $Q=(Q_0, Q_1, s, t)$ and denote by
$A:=kQ/\langle Q_2\rangle$ the radical square zero algebra of $Q$ over a field $k$.

\section{General theory for radical square zero algebras}\label{section3}
Let $k$ be a field. Let $A=kQ/\langle Q_2\rangle$ be the radical square zero algebra of a finite quiver $Q$ over $k$. Set $r=kQ_1$ and $E=kQ_0$.  It follows from Lemma 2.1 in \cite{Cib3} that 
%\begin{lemma}[Lemma 2.1 \cite{Cib3}]\label{lemma2.1}
%Let $k$ be a field. Let $r=kQ_1$ be the Jacobson radical of $A$ and  $E=kQ_0$.
 $r^{\otimes_E n}$ has a basis given by $Q_n$, the set of paths of length $n$. More precisely, 
for $p_i\in r$, we have that  $p_1\otimes_E \cdots \otimes_E p_n\neq 0\in r^{\otimes_E n}$ if and only if
  for $1\leq i \leq n-1$, $s(p_i)=t(p_{i+1}),$ that is, $p_1p_2\cdots p_n$ is a path of length $n$ in $kQ$.

Let $X$ and $Y$ be two sets of paths in $Q$. The set of {\it parallel paths} between $X$ and $Y$ is defined as
$$X//Y:=\{(\gamma, \gamma')\in X\times Y \ | \ \mbox{such that $s(\gamma)=s(\gamma')$ and $t(\gamma)=t(\gamma')$}  \}.$$
For instance, $Q_n//Q_0$ is the set of oriented cycles of length $n$. We  denote by $k(B)$, the vector space spanned by a given set $B$. Then we have the following result. 
\begin{lemma}[Lemma 2.2 \cite{Cib3}]\label{lemma2.2}
  For the algebra $A:=kQ/\langle Q_2\rangle$, the vector space $\Hom_{E-E}(r^{\otimes_E n}, A)$ is isomorphic to $k(Q_n//Q_0)\oplus k(Q_n//Q_1)$.
\end{lemma}
We denote by $\Omega_{nc}^i(A)$ the image of the differential $d_i$ in the projective resolution $\calR(A)$ (cf. (\ref{equ-proj-res}) above). In particular, $\Omega_{nc}^0(A)=A$. We remark that $\Omega_{nc}^i(A)$ is the bimodule of $i$-th noncommutative differential forms (relative to $E$) of $A$ (cf. \cite{CuQu}). 

\begin{lemma}\label{lemma-new34}
For any $p\in\Z_{\geq 0}$, there exists a natural isomorphism $\alpha: \Omega_{nc}^p(A)\rightarrow A\otimes_E  r^{\otimes_E p}$ of $A$-$A$-bimodules, where the right $A$-module structure on $A\otimes_E r^{\otimes_E p}$ is given by $$(a_0\otimes a_1\otimes \cdots\otimes a_p)\blacktriangleleft x=(\id^{\otimes p-1}\otimes \pi)\circ d_p(a_0\otimes a_1\otimes \cdots\otimes a_p\otimes x).$$ Here $\pi: A=E\oplus r\rightarrow r$ is the natural projection and for simplicity, we write $\otimes_E$ as $\otimes$.
\end{lemma}
\pf The map $\alpha$ is defined as the composition $$\Omega_{nc}^p(A) \hookrightarrow A\otimes r^{\otimes p-1}\otimes A \xrightarrow{\id^{\otimes p}\otimes \pi}A\otimes r^{\otimes p}.$$ Its inverse  $\alpha^{-1}: A\otimes r^{\otimes p}\rightarrow \Omega_{nc}^p(A)$ is given by $\alpha^{-1}(x)=(-1)^pd_p(x\otimes 1)$. By a direct computation, $\alpha$ is an isomorphism of $A$-$A$-bimodules.
\epf
\begin{lemma}\label{lemma-new}
 For any $p\in\Z_{\geq 0}$,  $\Omega_{nc}^p(A)$ has a basis given by $Q_{p}\cup Q_{p+1}$. As a consequence, $\Hom_{E-E}(r^{\otimes_E m}, \Omega_{nc}^p(A))$ is isomorphic to $k(Q_m//Q_{p})\oplus k(Q_m// Q_{p+1})$.
\end{lemma}
\pf This follows from  \cite[lemma2.1]{Cib3} and Lemma \ref{lemma-new34} since $A\otimes r^{\otimes p}$ is naturally isomorphic to $r^{\otimes p+1}\oplus r^{\otimes p}$. 
\epf 
\begin{rem}
In particular, when $p=0$ we recover the result in Lemma \ref{lemma2.2}. Lemma \ref{lemma-new} plays a crucial role in the present paper since it provides a combinatorial description of the Hochschild cochain complex $C^*(r, \Omega^p_{nc}(A))$ (cf. (\ref{equ-proj-res1}) above) for any $p\in\Z_{\geq 0}$. 
\end{rem}

 Let $D_{m, p}: k(Q_m//Q_p)\rightarrow k(Q_{m+1}//Q_{p+1})$ be the $k$-linear map defined by \begin{equation}\label{equ-defnD}
  D_{m, p}(\gamma_m, \gamma'_p):=\sum_{a\in Q_1}(a\gamma_m, a\gamma'_p)+(-1)^{p+m+1} \sum_{a\in Q_1} (\gamma_m a, \gamma'_pa).
\end{equation} 
Then  we have the following result. 
%Then we have the following Proposition.
\begin{prop}\label{prop-new}
  We have a commutative diagram,  where the vertical maps are given by the natural isomorphisms in Lemma \ref{lemma-new},
  \begin{equation}\label{1}
    \xymatrix@C=4pc{
    \Hom_{E-E}(r^{\otimes_E m}, \Omega_{nc}^p(A))\ar[d]^{\cong} \ar[r]^-{\delta} & \Hom_{E-E}(r^{\otimes_E m+1}, \Omega_{nc}^p(A))\ar[d]^{\cong}\\
    k(Q_m//Q_p)\oplus k(Q_m//Q_{p+1}) \ar[r]^-{\left( \begin{smallmatrix} 0 & D_{m,p} \\0 &0 \end{smallmatrix}\right)}& k(Q_{m+1}//Q_p)\oplus k(Q_{m+1}//Q_{p+1}).
    }
  \end{equation}
\end{prop}
\pf
Let $(\gamma_m, \gamma'_{p+1})\in Q_{m}//Q_{p+1}$. By the vertical isomorphism,  $(\gamma_m, \gamma'_{p+1})$ corresponds to the element
$\eta_{(\gamma_m, \gamma'_{p+1})}\in\Hom_{E-E}(r^{\otimes_E m}, \Omega_{nc}^p(A))$ defined as follows: For any $\alpha_m\in Q_m$,
\begin{equation*}
\eta_{(\gamma_m, \gamma_{p+1}')}(\alpha_m):=
\begin{cases}
\gamma_{p+1}'  & \mbox{if $\alpha_m=\gamma_m$;}\\
0  & \mbox{otherwise}
\end{cases}
\end{equation*}
Then from Remark \ref{rem-twoterms}, it follows that $\delta(\eta_{(\gamma_m, \gamma_{p+1}')})=0.$  Similarly, let $(\gamma_m, \gamma_{p}')\in Q_m//Q_p$. The vertical isomorphism sends $(\gamma_m, \gamma_{p}')$ to the element $\eta_{(\gamma_m, \gamma_{p}')}\in\Hom_{E-E}(r^{\otimes_E m}, \Omega_{nc}^p(A))$ defined by: 
\begin{equation*}
\eta_{(\gamma_m, \gamma_{p}')}(\alpha_m):=
\begin{cases}
e\gamma_{p}'  & \mbox{if $\alpha_m=\gamma_m$;}\\
0  & \mbox{otherwise.}
\end{cases}
\end{equation*}
From Remark \ref{rem-twoterms} again,  it follows that for any $\alpha_{m+1}\in Q_{m+1}$,
\begin{equation*}
\delta(\eta_{(\gamma_m, \gamma_{p}')})(\alpha_{m+1})=
\begin{cases}
a\gamma_p'  &  \mbox{if $\alpha_{m+1}=a\gamma_m$ for some $a\in Q_1$;}\\
(-1)^{p+m+1} \gamma_p'a &  \mbox{if $\alpha_{m+1}=\gamma_ma$ for some $a\in Q_1$;}\\
0 & \mbox{otherwise}.
\end{cases}
\end{equation*}
So the differential $\delta$ coincides with the map $D_{m, p}$. This proves the proposition. 
\epf

\begin{rem}\label{rem-coho}
  From Proposition \ref{prop-new}, it follows that
  \begin{equation}\label{equ-isom}
    \HH^m(A, \Omega_{nc}^p(A))\cong \frac{k(Q_m//Q_{p+1})}{\Imm(D_{m-1, p})}\oplus\Ker(D_{m,p}).
  \end{equation}
 Recall that we have  connecting homomorphisms  (cf. \cite[Section 3.1]{Wang})  $$\theta_{m, p}:\HH^m(A, \Omega_{nc}^p(A))\rightarrow \HH^{m+1}(A, \Omega_{nc}^{p+1}(A))$$  defined as $$ \theta_{m,p}(f)(a_{1, m+1})=f(a_{1, m})\otimes_E a_{m+1}.$$ Here we identify $\Omega_{nc}^p(A)$ with $A\otimes_E r^{\otimes_E p}$ by Lemma \ref{lemma-new34}.
\end{rem}

\begin{lemma}\label{lemma-limitnew}
  Under the isomorphism (\ref{equ-isom}) above, we have the following commutative diagram.
  \begin{equation}\label{diagram3}
    \xymatrix@C=5pc{
    \HH^m(A, \Omega_{nc}^p(A)) \ar[r]^-{\theta_{m,p}}\ar[d]^{\cong} & \HH^{m+1}(A, \Omega_{nc}^{p+1}(A))\ar[d]^{\cong}\\
    \frac{k(Q_m//Q_{p+1})}{\Imm(D_{m-1, p})}\oplus\Ker(D_{m,p})\ar[r]^-{\left( \begin{smallmatrix} E_{m, p+1} & 0\\ 0 &E_{m, p} \end{smallmatrix}\right)} & \frac{k(Q_{m+1}//Q_{p+2})}{\Imm(D_{m, p+1})}\oplus\Ker(D_{m+1,p+1})
    }
  \end{equation}
  where $E_{m, p+1}: \frac{k(Q_m//Q_{p+1})}{\Imm(D_{m-1, p})}\rightarrow \frac{k(Q_{m+1}//Q_{p+2})}{\Imm(D_{m, p+1})}$ is defined as,
  \begin{equation*}
    E_{m, p+1}(\gamma_m, \gamma_{p+1}')=-\sum_{a\in Q_1} (\gamma_ma,\gamma_{p+1}'a).
  \end{equation*}
 % and $$F_{m,p}: \Ker(D_{m,p})\rightarrow\Ker(D_{m+1,p+1})$$
 % is defined,
 % \begin{equation*}
 %   F_{m,p}(\gamma_m, \gamma_p')=-\sum_{a\in Q_1} (\gamma_ma, \gamma_p'a).
%  \end{equation*}
\end{lemma}
\pf First, let us prove that $E_{m, p+1}$ is well-defined on $ \frac{k(Q_m//Q_{p+1})}{\Imm(D_{m-1, p})}$ and $E_{m, p}$ is well-defined on $\Ker(D_{m, p})$. Let $(\gamma_{m-1}, \gamma_p')\in Q_{m-1}//Q_p$. We need to prove that $$E_{m, p+1}(D_{m-1, p}((\gamma_{m-1}, \gamma_p'))\in \Imm(D_{m, p+1}).$$ Indeed, we have
\begin{equation*}
\begin{split}
E_{m, p+1}(D_{m-1, p}((\gamma_{m-1}, \gamma_p'))&=\sum_{a\in Q_1} E_{m, p+1}(a\gamma_{m-1}, a\gamma_p')+(-1)^{p+m}
\sum_{a\in Q_1} E_{m, p+1}(\gamma_{m-1}a, \gamma_p'a)\\
&=\sum_{a, b\in Q_1} -(a\gamma_{m-1}b, a\gamma_p'b)+(-1)^{p+m+1}\sum_{a, b\in Q_1}
(\gamma_{m-1}ab, \gamma_p'ab)\\
&=\sum_{a\in Q_1} D_{m, p+1}(\gamma_{m-1}a, \gamma_{p}'a).
\end{split}
\end{equation*}
Similarly, we can also prove that $E_{m, p}$ is well-defined on $\Ker(D_{m, p})$.

Let $(\gamma_m, \gamma_{p+1}')\in Q_m//Q_{p+1}$.
Then by the vertical isomorphism, $(\gamma_m, \gamma_{p+1}')$ is sent to
an element in $\HH^m(A, \Omega_{nc}^p(A))$ represented by  $\eta_{(\gamma_m, \gamma_{p+1}')}\in\Hom_{E-E}(r^{\otimes_E m}, \Omega_{nc}^p(A)),$
where $\eta_{(\gamma_m, \gamma_{p+1}')}$ has been defined in the proof of Proposition \ref{prop-new}.
From (\ref{equ-isom}), it follows that
\begin{equation*}
\theta_{m, p}(\eta_{(\gamma_m, \gamma_{p+1}')})(a_{m+1})=
\begin{cases}
-\gamma_{p+1}'a & \mbox{if $a_{m+1}=\gamma_m a$ for some $a\in Q_1$;}\\
0 &  \mbox{otherwise}.
\end{cases}
\end{equation*}
Similarly, let $z:=\sum_{(\gamma_m, \gamma_p')\in (Q_m//Q_p)} x_{(\gamma_m, \gamma_p')} (\gamma_m, \gamma_p')$ be an element in $\Ker(D_{m, p})$. Then $z$ is sent to an element in $\HH^m(A, \Omega_{nc}^p(A))$ represented by
$$\eta_z:=\sum_{(\gamma_m, \gamma_p')\in (Q_m//Q_p)} x_{(\gamma_m, \gamma_p')} \eta_{(\gamma_m, \gamma_p')}\in\Hom_{E-E}(r^{\otimes_E m}, \Omega_{nc}^p(A)). $$
%$$\theta_{m, p}(
This implies
\begin{equation*}
\theta_{m, p}(\eta_z)(a_{m+1})=
\begin{cases}
-\gamma_pa &  \mbox{if $a_{m+1}=\gamma_ma$ for some $a\in Q_1$;}\\
0 & \mbox{otherwise.}
\end{cases}
\end{equation*}
This completes the proof.
\epf

\begin{prop}
  For $m\in\Z$, we have
  \begin{equation}\label{prop-isom}
    \HH_{\sg}^m(A, A)\cong \lim_{\substack{\longrightarrow\\p\in\Z_{> 0}\\ m+p>0}}\frac{k(Q_{m+p}//Q_{p+1})}{\Imm(D_{m+p-1, p})}\oplus \lim_{\substack{\longrightarrow\\p\in\Z_{> 0}\\ m+p>0}} \Ker(D_{m+p, p}),
  \end{equation}
  where the colimits are taken along maps $E_{m, p}$. 
\end{prop}
\pf This follows from  Lemma \ref{lemma-limitnew} and \cite[Section 3.2]{Wang}.\epf

\begin{cor}
   If the quiver $Q$ has no oriented cycles, then $\HH_{\sg}^m(A, A)=0$ for any $m\in\Z$.
\end{cor}
\pf If $Q$ has no oriented cycles, then for $p\gg0$, we have $k(Q_{m+p}//Q_{p+1})=0.$ Thus the right hand side of the isomorphism
(\ref{prop-isom}) vanishes. This yields $\HH^m_{\sg}(A, A)=0$.
\epf
\begin{rem}
  In general, the map $D_{m, p}$ is not injective for $m, p\in\Z_{>0}$. In the following,  we will study  in which case $D_{m, p}$ is injective.
\end{rem}
\begin{defn}\label{defn-crown}
  Let $c\in \Z_{\geq 0}$. A $c$-crown is a quiver with $c$-vertices
  cyclically labelled by the cyclic group of order $c$, and $c$ arrows
  $a_0, \cdots, a_{c-1}$ such that $s(a_i)=i$ and $t(a_i)=i+1$. For
  instance, a $1$-crown is a loop, and a $2$-crown is the two-way quiver.
\end{defn}
Now let us consider a finite and connected quiver $Q=(Q_0, Q_1, s, t)$ without sources or sinks (cf. e.g.
\cite{Mich}), that is,  for any vertex $e\in Q_0$, there
exist arrows $p, q\in Q_1$ such that $s(p)=t(q)=e$.
For instance, $c$-crowns and the following quivers have no sources or sinks.

$$
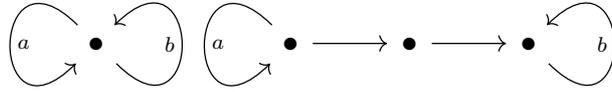
\begin{tikzcd}
% \bullet  \ar[,loop, out=315, in=45, distance=3em]{}{a} & &
 \bullet \ar[%
    ,loop % tells tikz-cd to do a loop
    ,out=135 % start at angle 123?
    ,in=225 % stop at angle 57?
    ,distance=4em % biggest distance of arrow to node. Yarou can use pt or cm as well.
    ]{}{a} \ar[
    , loop
    ,out=315
    ,in=45
    ,distance=4em]{}{b}
 &  & \bullet \arrow[ loop, out=135, in=225, distance=4em]{}{a} \ar{r}&\bullet \arrow{r} & \bullet \arrow[loop, out=315,in=45, distance=4em ]{}{b}
\end{tikzcd}$$

%\begin{lemma}\label{lemma-inter}
%  Let $Q=(Q_0, Q_1, s, t)$ be a finite and connected quiver without sources or sinks.
%  For any two vertices $e_0, e_1\in Q_0$ (possibly $e_0=e_1$), there exists a path
%  $\beta$ such that $s(\beta)=e_0$ and $t(\beta)=e_1$.
%\end{lemma}
%\pf

%For such quiver, we observe that $D_{m,p}$ is injective for almost $m, p\in \Z_{>0}$.
%Namely, we have the following Proposition.
\begin{prop}\label{prop-com}
  Let $k$ be a field. Let $Q=(Q_0, Q_1, s, t)$ be a finite and connected quiver without sources or sinks. Suppose that $Q$ is not a crown. Then we have the following two cases for $m,p\in \Z_{>0}$:
  \begin{enumerate}
      \item If $m=p$, then we have
    $\Ker(D_{m, m})$ is a one-dimensional $k$-vector space with a basis $\sum_{\gamma_m\in Q_{m}} (\gamma_m, \gamma_m)$.
  \item If $m\neq p$, then the map $D_{m, p}: k(Q_m//Q_p)\rightarrow k(Q_{m+1}//Q_{p+1})$ is injective.
  \end{enumerate}
\end{prop}
\pf This proof is analogous to the one of \cite[Theorem 3.1]{Cib3}. Suppose that $$x=\sum_{(\gamma_m, \beta_p)\in Q_m//Q_p}x_{(\gamma_m, \beta_p)}
(\gamma_m, \beta_p)$$
is an element in $\Ker(D_{m, p})$, where $x_{(\gamma_m, \beta_p)}\in k$ and those $(\gamma_m, \beta_p)$ in the above sum are distinct. Let us fix an element
$(\gamma_m, \beta_p)$ such that $x_{(\gamma_m, \beta_p)}\neq 0$.
We consider the contribution of $x_{(\gamma_m, \beta_p)}(\gamma_m, \beta_p)$
in $D_{m, p}(x)$, which is
\begin{equation}\label{equ-contribution}
  x_{(\gamma_m, \beta_p)}\left(\sum_{a\in Q_1s(\gamma_m)}(a\gamma_m, a\beta_p)+(-1)^{p+m+1} \sum_{a\in t(\gamma_m)Q_1} (\gamma_m a, \beta_pa)\right).
\end{equation}
Then from (\ref{equ-contribution}), we have the following observation.
\begin{obs}\label{obs}
  If $x_{(\gamma_m, \beta_p)}\neq 0$, then for $1\leq i\leq \min\{m, p\}$, $f_i(\gamma_m)=f_i(\beta_p)$ and $l_i(\gamma_m)=l_i(\beta_p),$ where $f_i(\gamma)$ and $l_i(\gamma)$ denote the first $i$ arrows and last $i$ arrows of a path $\gamma$ respectively.
\end{obs}
Suppose $m=p$.
%Let $$x=\sum_{(\gamma_m, \beta_m)\in k(Q_m//Q_m)}x_{(\gamma_m, \beta_m)}
%(\gamma_m, \beta_m)$$
%be an element in $\Ker(D_{m, m}).$ Fix an element
%$(\gamma_m, \beta_m)$ such that $x_{(\gamma_m, \beta_m)}\neq 0$.
From Observation \ref{obs} above, it follows that
$x_{(\gamma_m, \beta_m)}\neq 0$ if and only if $\gamma_m=\beta_m$. Hence we can write $x\in \Ker(D_{m,m})$ as $$x=\sum_{(\gamma_m, \gamma_m)\in Q_m//Q_p}x_{(\gamma_m, \gamma_m)} (\gamma_m, \gamma_m).$$
We claim that, for any $(\gamma_m, \gamma_m), (\gamma_m', \gamma_m')\in
Q_{m}//Q_m$, $x_{(\gamma_m, \gamma_m)}=x_{(\gamma_m', \gamma_m')}.$  Indeed, we define an equivalence relation on $Q_m$ as follows. For $$\gamma_m:=a_1\cdots a_m, \gamma_m':=b_1\cdots b_m\in Q_m,$$ we say that
$\gamma_m\sim \gamma_m'$ if we have $a_{i}=b_{i+1}$ for $1\geq i \geq m-1,$ or $b_{i}=a_{i+1}$ for $1\geq i\geq m-1$. Then we can extend $\sim$ to an equivalence relation $\sim_{c}$ on $Q_m$. Observe that if $\gamma_m\sim_c \gamma_m'$, then $x_{(\gamma_m, \gamma_m)}= x_{(\gamma_m', \gamma_m')}.$  So it is sufficient to show that for any $\gamma_m, \gamma_m'\in Q_m$, $\gamma_m\sim_c \gamma_m'.$ Since $Q$ is a finite and connected quiver without sources or sinks, we have the following two observations.
\begin{obs}\label{obs-3.14}
  If $\gamma_m$ intersects with $\gamma_m'$ (that is,
  there exist a vertex $e\in Q_0$ such that $e\in \gamma_m$ and $e\in \gamma_m'$), then $\gamma_m\sim_c \gamma_m'.$
\end{obs}
\begin{obs}\label{obs-3.15}
  Given any path $\beta$ of length smaller than $m$, we can extend $\beta$ to
  a path of length $m$.
\end{obs}
Now let us prove $\gamma_m\sim_c\gamma_m'$ for any $\gamma_m, \gamma_m'\in Q_m$ by
the two observations above. Since $Q$ is finite and connected,
there exists a non-oriented path connecting $\gamma_m$ and $\gamma_m'$.
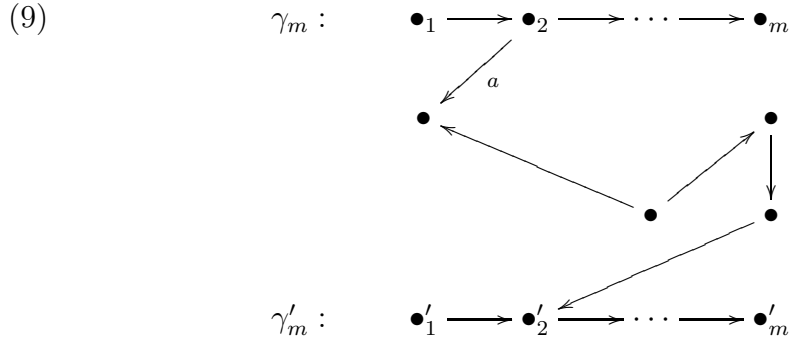
\begin{equation} \label{dia3.14}
  \xymatrix{
  \gamma_m:&  \bullet_1 \ar[r]  & \bullet_2\ar[dl]^{a}\ar[r]  &\cdots\ar[r] & \bullet_m \\
&\bullet&  &&\bullet\ar[d]&&\\
       &  &&  \bullet\ar[ull]\ar[ur] &\bullet\ar[dll]&&\\
\gamma_m': & \bullet_1'\ar[r]  & \bullet_2'\ar[r]  &\cdots\ar[r]  & \bullet_m'\\
  }
\end{equation}
Let us use induction on the length $l(\beta)$ of the (non-oriented) path $\beta$ connecting $\gamma_m$ and $\gamma_m'$. It is clear for $l(\beta)=0$ by Observation \ref{obs-3.14}. Assume that $\gamma_m\sim_c\gamma_m'$ for any $\gamma_m, \gamma_m'\in Q_m$ such that there exists a path of length smaller than $(l-1)$ between them. Now by Observation \ref{obs-3.15}, the arrow $a$ in Diagram (\ref{dia3.14}) can be extended to a path $\gamma_m^a$ of length $m$. From Observation (\ref{obs-3.14}) again it follows that $\gamma_m\sim_c \gamma_m^a.$ Note that the length of the (non-oriented) path connecting $\gamma_m'$ and $\gamma_m^a$ is $(l-1)$, thus by induction hypothesis, we have $\gamma_m^a\sim_c \gamma_m'.$ Since $\sim_c$ is transitive, we have $\gamma_m\sim_c \gamma_m'$ for any $\gamma_m, \gamma_m'\in Q_m$.
%Since $Q$ is a finite and connected quiver without sources or sinks,
%from Lemma \ref{lemma-inter} it follows that
%there exist a path $\beta$ such that
%$$s(\beta)=t(\gamma_m')$$ and $$t(\beta)=s(\gamma_m).$$
%Hence along the path $\beta$, we can induce
%$\gamma_m\sim_c \gamma_m'$.
Therefore we have $x\in k\left(\sum_{\gamma_m\in Q_{m}} (\gamma_m, \gamma_m)\right)$ if $x\in \Ker(D_{m,m}).$ On the other hand, for any $\lambda\in k$, we observe that $D_{m,m}(\lambda \sum_{\gamma_m\in Q_m}(\gamma_m, \gamma_m))=0.$ Indeed,
\begin{equation*}
  \begin{split}
    D_{m, m}(\lambda \sum_{\gamma_m\in Q_m}(\gamma_m, \gamma_m))&=\lambda \sum_{a\in Q_1t(\gamma_m)}\sum_{\gamma_m\in
    Q_m}(a\gamma_m, a\gamma_m)-\lambda \sum_{a\in s(\gamma_m)Q_1}\sum_{\gamma_m\in Q_m}
    (\gamma_ma, \gamma_ma)\\
    &=\lambda \sum_{\gamma_{m+1}\in Q_{m+1}}(\gamma_{m+1}, \gamma_{m+1})-\lambda
    \sum_{\gamma_{m+1}\in Q_{m+1}}(\gamma_{m+1}, \gamma_{m+1})\\
    &=0.
  \end{split}
\end{equation*}
Therefore $\Ker(D_{m, m})=k\left(\sum_{\gamma_m\in Q_{m}} (\gamma_m, \gamma_m)\right).$

Suppose $m\neq p$. Without loss of generality, we may assume that $m<p$. Then there will be two cases.
\begin{enumerate}
\item If $p\geq 2m$, let us write $\gamma_m=a_1a_2\cdots a_m$ (possibly $a_i=a_j$ for some $i, j\in\{1, 2, \cdots, m\}$). If $x_{(\gamma_m, \beta_p)}\neq 0$, then from Observation \ref{obs} it follows that
 $$\beta_p=a_1\cdots a_mb_1b_2\cdots b_{p-2m} a_1\cdots a_m$$ is an oriented cycle (possibly with self-intersection).
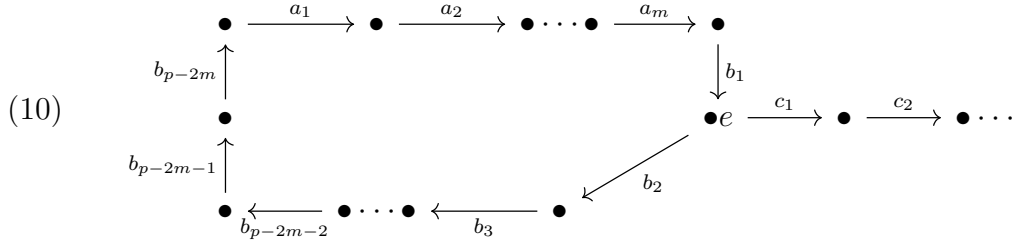
\begin{equation}\label{equ-qui}
\begin{tikzcd}
\bullet \ar{r}{a_1}& \bullet \ar{r}{a_2} & \bullet \cdots\bullet\ar{r}{a_m} & \bullet \ar{d}{b_1} &\\
\bullet \ar{u}{b_{p-2m}} & &&\bullet e\ar{dl}{b_2}\ar{r}{c_1} &\bullet\ar{r}{c_2}& \bullet \cdots& \\
\bullet\ar{u}{b_{p-2m-1}}&\bullet \cdots\bullet\ar{l}{b_{p-2m-2}} & \bullet\ar{l}{b_3}&
\end{tikzcd}
\end{equation}
Since $Q$ is not a crown, there exists a vertex, say $e$,  such that there is an another arrow $c_1$ starting at $e$
and $c_1\neq b_2$ (or there is an another arrow $c_1$ ending at $e$ and $c_1\neq b_1$, respectively)  (cf. Diagram (\ref{equ-qui})). Let us just consider the first case, namely where there is an arrow $c_1$ starting at $e$ and $c_1\neq b_2$, the other case being analogous. From (\ref{equ-contribution}), we obtain that $$(-1)^{p+m+1}x_{(\gamma_m, \beta_p)}(a_1\cdots a_mb_1, a_1\cdots a_mb_1\cdots b_{p-2m}a_1\cdots a_m b_1)$$ is a nonzero summand of $D_{m,p}(x)$.
Since $D_{m, p}(x)=0$, we get that
$x_{(\gamma_m', \beta_p')}\neq 0,$
where $(\gamma_m',\beta_p'):=(a_2\cdots a_mb_1, a_2\cdots a_mb_1\cdots b_{p-2m}a_1\cdots a_m b_1).$
Thus it follows that
$$(-1)^{p+m+1}x_{(\gamma_m', \beta_p')}(a_2\cdots a_mb_1c_1, a_2\cdots a_mb_1\cdots b_{p-2m}a_1\cdots a_m b_1c_1)$$ is a nonzero summand of $D_{m,p}(x).$ By this induction process, we have that $$x_{(c_1c_2\cdots c_{m-1}c_m, b_2\cdots b_{p-2m}a_1\cdots a_mb_1c_1\cdots c_{m-1}c_m)}\neq 0.$$ Since $c_1\neq b_2$, from Observation \ref{obs} it follows that
$$x_{(c_1c_2\cdots c_{m-1}c_m, b_2\cdots b_{p-2m}a_1\cdots a_mb_1c_1\cdots c_{m-1}c_m)}=0.$$
Contradiction! Therefore, we have $x=0$.
For the remaining cases, we can apply an analogous induction process to
argue $x=0$.
%\item If $m< p< 2m$, the proof is similar to the one of Case 1. From Observation \ref{obs}.
%it follows that $\gamma_m$ is an oriented cycle (possibly with self-intersection).
%\begin{equation}\label{equ-qui2}
%\begin{tikzcd}
%e\bullet \ar{r}{a_1} & \bullet\ar{r}{a_2} & \bullet \ar{d}{a_3} &\\
%\bullet \ar{u}{a_{r}} &\bullet\cdots\bullet \ar{l}{a_{r-1}} &\bullet\ar{r}{c_1}\ar{l}{a_4}& \bullet \cdots& \\
%\end{tikzcd}
%\end{equation}
%Suppose $s(\gamma_m)=t(\gamma_m)=e$ in Diagram (\ref{equ-qui2}).
%Since $Q$ is not a crown, there exists an arrow $c_1$ such that $c_1\neq a_4$, as we have shown in Diagram (\ref{equ-qui2}).
%Then by the induction process similar to the one in Case 1,
%we can show that $x=0$.

\end{enumerate}
Therefore, we have proved that $\Ker(D_{m,p})=0$ for $m\neq p$. The proof is complete.
\epf

%So from Proposition \ref{prop-com}, for a quiver $Q$ without sources or sinks, the Tate-Hochschild cohomology $\HH_{\sg}^m(A, A)$ of its radical
%square zero algebra $A$ has a
%much combinatorial description.
\begin{thm}\label{prop-singular}
Let $k$ be a field. Let $Q$ be a finite and connected quiver without sources or sinks and $A$ be its radical square
zero $k$-algebra. Suppose that $Q$ is not a crown. Then we have the following:
\begin{enumerate}
\item For $m, p\in\Z_{>0}$,
\begin{equation*}
\HH^m(A, \Omega_{nc}^p(A))=
\begin{cases}
\frac{k(Q_m//Q_{p+1})}{\Imm(D_{m-1, p})} & \mbox{if $m\neq p$},\\
\frac{k(Q_m//Q_{m+1})}{\Imm(D_{m-1,m})}\oplus \Ker(D_{m, m}) & \mbox{if $m=p$}.
\end{cases}
\end{equation*}
\item The connecting homomorphism $\theta_{m, p}: \HH^m(A, \Omega_{nc}^p(A))\rightarrow \HH^{m+1}(A, \Omega_{nc}^{p+1}(A))$ is injective for any $m, p\in \Z_{>0}$.
\item The Tate-Hochschild cohomolgy ring $\HH^*_{\sg}(A, A)$ has a filtration
$$0\subset \HH^*(A, A)\subset \cdots\subset \HH^{*+p}(A, \Omega_{nc}^p(A))\subset \HH^{*+p+1}(A, \Omega_{nc}^{p+1}(A))\subset \cdots$$
Moreover, this filtration respects the Gerstenhaber algebra structure, that is, for $m, n\in \Z$,
$$\HH^{m+p}(A, \Omega_{nc}^p(A))\cup\HH^{n+q}(A, \Omega_{nc}^q(A))\subset \HH^{m+p+n+q}(A, \Omega_{nc}^{p+q}(A)),$$
$$[\HH^{m+p}(A, \Omega_{nc}^p(A)), \HH^{n+q}(A, \Omega_{nc}^q(A))]\subset \HH^{m+p+n+q-1}(A, \Omega_{nc}^{p+q}(A)).$$
\end{enumerate}
\end{thm}
\pf The first assertion follows from Proposition \ref{prop-com}. Thus  it is sufficient to verify that $\theta_{m,p}$ is injective for $m,p\in\Z_{>0}$.
Recall that we have the long exact sequence
\begin{equation*}
\xymatrix{
\cdots\ar[r] &\HH^m(A, A\otimes r^{\otimes p}\otimes A)\ar[r]^-{d} &\HH^m(A, \Omega_{nc}^p(A))\ar[r]^-{\theta_{m, p}} & \HH^{m+1}(A, \Omega_{nc}^{p+1}(A))\ar[r]& \cdots
}
\end{equation*}
where $d$ is induced by
the differential of the resolution $\calR(A)$ (cf. Section \ref{section2}).
Hence $\theta_{m, p}$ is injective if and only if
$d=0$. Now let us show that $d=0$ for any $m, p\in\Z_{>0}$.
Note that $\Hom_{E-E}(r^{\otimes_E m}, A\otimes_Er^{\otimes_E p}\otimes_E A)$ has a decomposition
with respect to the decomposition $A\cong E\oplus r$. Namely,
\begin{equation*}
\begin{split}
\Hom_{E-E}(r^{\otimes_E m}, A\otimes_Er^{\otimes_E p}\otimes_E A)\cong&\Hom_{E-E}(r^{\otimes_E m}, r\otimes_Er^{\otimes_E p}\otimes_E r)\\
&\bigoplus\Hom_{E-E}(r^{\otimes_E m}, E\otimes_Er^{\otimes_E p}\otimes_E r)\\
&\bigoplus\Hom_{E-E}(r^{\otimes_E m}, r\otimes_Er^{\otimes_E p}\otimes_E E)\\
&\bigoplus\Hom_{E-E}(r^{\otimes_E m}, E\otimes_Er^{\otimes_E p}\otimes_E E),
\end{split}
\end{equation*}
hence $\Hom_{E-E}(r^{\otimes_E m}, A\otimes_Er^{\otimes_E p}\otimes_E A)$ has a basis $$S_{m,p}:=(Q_m// Q_{p+2})\cup (Q_m//eQ_{p+1})\cup (Q_m//Q_{p+1}e)\cup (Q_m//Q_{p}),$$
where we use the word $e$ to distinguish $(Q_m//eQ_{p+1})$ and $(Q_m//Q_{p+1}e)$, more precisely,
$(Q_m//eQ_{p+1})$ is the basis corresponding to $\Hom_{E-E}(r^{\otimes_E m}, E\otimes_Er^{\otimes_E p}\otimes_E r)$ and
$(Q_m//Q_{p+1}e)$ is the one corresponding to $\Hom_{E-E}(r^{\otimes_E m}, r\otimes_Er^{\otimes_E p}\otimes_E E)$.
Moreover, we can write down the differential.
\begin{equation}\label{diffen-equ}
\xymatrix{
\Hom_{E-E}(r^{\otimes_E m}, A\otimes_Er^{\otimes_E p}\otimes_E A)\ar[d]^{\cong}\ar[r]^-{\delta}&\Hom_{E-E}(r^{\otimes_E m+1}, A\otimes_Er^{\otimes_E p}\otimes_E A)\ar[d]^{\cong}\\
k(S_{m, p}) \ar[r]_-{\left( \begin{smallmatrix} 0 & D_{m, p+1}&D_{m,p+1}' &0\\0& 0 & 0 & E_{m, p}\\0& 0 & 0 & E_{m, p}'\\
0 & 0 & 0 & 0 \end{smallmatrix}\right)} & k(S_{m+1, p})
 }\end{equation}
where
\begin{equation*}
\begin{split}
D_{m, p+1}(\gamma_m, e\beta_{p+1}):=&\sum_{a\in Q_1}(a\gamma_m, a\beta_{p+1});\\
D_{m, p+1}'(\gamma_m, \beta_{p+1}e):=&\sum_{a\in Q_1}(-1)^{m+1}(\gamma_ma, \beta_{p+1}a);\\
E_{m, p}(\gamma_m, \beta_{p}):=& \sum_{a\in Q_1}(a\gamma_m, a\beta_{p+1});\\
E_{m, p}'(\gamma_m,\beta_{p}):=& \sum_{a\in Q_1}(-1)^{m+1}(\gamma_ma, \beta_{p+1}a).
\end{split}
\end{equation*}
Now suppose $$x\in \HH^m(A, A\otimes r^{\otimes p}\otimes A)$$
is a nonzero element.
Then from Diagram (\ref{diffen-equ}), we can write $x$ as follows,
\begin{equation*}
\begin{split}
x=&\sum_{(\gamma_m, \beta_{p+2})\in Q_m//Q_{p+2}} x_{(\gamma_m, \beta_{p+2})}(\gamma_m, \beta_{p+2})+\sum_{(\gamma_m, \beta_{p+1})\in Q_m//Q_{p+1}}  x_{(\gamma_m, e\beta_{p+1})}(\gamma_m, e\beta_{p+1})
+\\
&\sum_{(\gamma_m, \beta_{p+1})\in Q_m//Q_{p+1}}  x_{(\gamma_m, \beta_{p+1}e)}(\gamma_m, \beta_{p+1}e).
\end{split}
\end{equation*}
Now $\delta(x)=0$ is equivalent to
\begin{equation}\label{equ-contri}
\sum_{\substack{ (\gamma_m, \beta_{p+1})\in (Q_m//Q_{p+1}) \\a\in Q_1}} \left( x_{(\gamma_m, e\beta_{p+1})}(a\gamma_m, a\beta_{p+1})
+(-1)^{p+m+1}x_{(\gamma_m, \beta_{p+1}e)}(\gamma_ma, \beta_{p+1}a)\right)=0.
\end{equation}
We have the following observation.
\begin{obs}\label{obs3}
Suppose that  $x_{(\gamma_m, e\beta_{p+1})}\neq 0.$
Then
$t(\gamma_m)=t(\beta_{p+1})$.
%there exists $(\gamma_m', \beta_{p+1}'e)$ such that
%$$x_{(\gamma'_m, e\beta'_{p+1})}\neq 0$$
Moreover, for those $a\in Q_1$ such that $a\gamma_{m}\neq 0$, we have
$$x_{(af_{m-1}(\gamma_{m}), af_{p}(\beta_{p+1}e)}=(-1)^{p+m}x_{(\gamma_m, e\beta_{p+1})},$$
where we denote by $f_{m-1}(\gamma_{m})$ the path formed by the first $m-1$ arrows in $\gamma_m$. Similarly, suppose that $x_{(\gamma'_m, \beta'_{p+1}e)}\neq 0.$ Then $s(\gamma'_m)=s(\beta'_{p+1}).$ Moreover, for those $a\in Q_1$ such that  $\gamma_m'a\neq 0$, we have $$x_{(l_{m-1}(\gamma_{m}')a, el_{p}(\beta_{p+1}')a)}=(-1)^{m+p}x_{(\gamma_m', \beta_{p+1}'e)},$$
where we denote by $l_{m-1}(\gamma'_{m})$ the path formed by the last $m-1$ arrows in $\gamma'_m$.
\end{obs}

Now let us compute
\begin{equation*}
\begin{split}
d(x)=\sum x_{(\gamma_m, e\beta_{p+1})}(\gamma_m, \beta_{p+1})
+\sum (-1)^p x_{(\gamma_m, \beta_{p+1}e)}(\gamma_m, \beta_{p+1})\in\HH^m(A, \Omega_{nc}^p(A)).
\end{split}
\end{equation*}
We claim that $d(x)=(-1)^{m}\sum_{(\gamma_m, e\beta_{p+1})} x_{(\gamma_m, e\beta_{p+1})}D_{m-1, p}(f_{m-1}(\gamma_{m}), f_p(\beta_{p+1})).$
Indeed, we have that
\begin{eqnarray*}
\lefteqn{
\sum_{(\gamma_m, e\beta_{p+1})} x_{(\gamma_m, e\beta_{p+1})}D_{m-1, p}(f_{m-1}(\gamma_{m}), f_p(\beta_{p+1}))}\\
&=&\sum_{\substack{ (\gamma_m,e \beta_{p+1}), \\a\in Q_1}} x_{(\gamma_m, e\beta_{p+1})}(af_{m-1}(\gamma_{m}), af_p(\beta_{p+1}))+\\
&&(-1)^{p+m}\sum_{\substack{ (\gamma_m,e \beta_{p+1}), \\a\in Q_1}}x_{(\gamma_m, e\beta_{p+1})}(f_{m-1}(\gamma_{m})a, f_p(\beta_{p+1})a)\\
&=&\sum_{\substack{ (\gamma_m,e \beta_{p+1}), \\a\in Q_1}}(-1)^{p+m} x_{(af_{m-1}(\gamma_{m}), af_p(\beta_{p+1})e)}(af_{m-1}(\gamma_{m}), af_p(\beta_{p+1}))+\\
&&(-1)^m\sum_{\substack{ (\gamma_m,e \beta_{p+1}), \\a\in Q_1}}x_{(f_{m-1}(\gamma_m)a,  ef_p(\beta_{p+1})a)}(f_{m-1}(\gamma_{m})a, f_p(\beta_{p+1})a)\\
&=&(-1)^md(x),
\end{eqnarray*}
where the second identity comes from Observation \ref{obs}.
So it follows that  $d(x)=0$ in $\HH^m(A, \Omega_{nc}^p(A))$. Therefore we have showed that $\theta_{m, p}$ is injective. Assertion (3) follows from (\ref{defn-bracket}) below. This proves the theorem.
\epf

As a corollary, we have the following result.
\begin{cor}\label{cor-cup}
Let $Q$ be a finite and connected quiver without sources or sinks and $A$ be its radical square
zero $k$-algebra over a field $k$. Suppose that $Q$ is not a crown. Then for any $m, n, p, q\in\Z_{>0}$ such that $m\neq p$ and $n\neq q$, we have
$$\HH^m(A, \Omega_{nc}^p(A))\cup\HH^n(A, \Omega_{nc}^q(A)=0.$$
In particular, for $m, n\in \Z$ and $mn\neq 0$, we have 
$$\HH^m_{\sg}(A, A)\cup \HH^n_{\sg}(A, A)=0.$$
\end{cor}
\pf From Theorem \ref{prop-singular}, it follows that $$\HH^m(A, \Omega_{nc}^p(A))=\frac{k(Q_m//Q_{p+1})}{\Imm(D_{m-1, p})}$$ for $m\neq p$. Note that we have the following identity on the level of chains:  For $m\neq p$ and $n\neq q$, $$k(Q_m//Q_{p+1})\cup k(Q_n//Q_{q+1})=0.$$ Hence on the level of cohomology groups, we have  $$\HH^m(A, \Omega_{nc}^p(A))\cup \HH^n(A, \Omega_{nc}^q(A))=0$$ for $m\neq p$ and $n\neq q$.  In particular, $\HH_{\sg}^m(A, A)\cup \HH^n_{\sg}(A, A)=0,$ for $mn\neq 0$.
\epf
%\begin{rem}\label{rem-conter}
%From Corollary \ref{cor-cup} above,  it follows that the Gerstenhaber brackets on $\HH^*(A, A)$ and $\HH_{\sg}^*(A, A)$ cannot come from a BV structure since
%the cup products vanish and the Lie brackets do not vanish in general (cf. Section \ref{two loops quiver}).
%\end{rem}

In the rest of this section, we will study the Gerstenhaber bracket $[\cdot,\cdot]$ (cf. \cite{Wang}) on the total complex $$\bigoplus_{m, p\in \Z_{\geq 0}}\Hom_{E-E}(r^{\otimes_E m}, \Omega_{nc}^p(A)).$$ For $m, p \in Z_{>0}$,  we denote $C^m(r, \Omega_{nc}^p(A)):=\Hom_{E-E}(r^{\otimes_E m}, A\otimes_E r^{\otimes_E p}).$ Recall that we identify $\Omega_{nc}^p(A)$ with $A\otimes_E r^{\otimes_E p}$ by Lemma \ref{lemma-new34}.  

Take two elements $f\in C^m(r, \Omega_{nc}^p(A))$ and $g\in C^n(r, \Omega_{nc}^q(A))$.
Denote
\begin{equation*}
f\bullet_i g:=
\begin{cases}
(f\otimes \id^{\otimes q})(\id^{\otimes i-1}\otimes \overline{g}\otimes \id^{\otimes m-i})&\mbox{if} \ 1\leq i\leq m, \\
(\id^{\otimes -i}\otimes \overline{g}\otimes \id^{\otimes p+i})(f\otimes \id^{\otimes n-1}) & \mbox{if} \ -p\leq i \leq -1,
\end{cases}
\end{equation*}
where $\otimes$ represents $\otimes_E$ and $\overline{g}=(\id^{\otimes p}\otimes \pi)\circ g$.
We also denote
\begin{equation*}
f\bullet g:=\sum_{i=1}^m(-1)^{(i-1)(q-n-1)}f\bullet_i g -\sum_{i=1}^p(-1)^{(i-m-p-1)(q-n-1)}f\bullet_{-i} g.
\end{equation*}
Then we define
\begin{equation}\label{defn-bracket}
  [f, g]:=f\bullet g-(-1)^{(m-p-1)(n-q-1)} g\bullet f.
\end{equation}
Note that $[f, g]\in C^{m+n-1}(r, \Omega_{nc}^{p+q}(A)).$ Then from  \cite{Wang}, it follows that $[\cdot, \cdot]$ defines a differential graded Lie algebra structure on the total complex
$$\bigoplus_{m\in \Z_{>0}, p\in \Z_{\geq 0}}C^m(r, \Omega_{nc}^p(A))$$
and thus $[\cdot, \cdot]$ defines a graded Lie algebra structure on the cohomology groups
$$\bigoplus_{m\in\Z_{>0}, p\in \Z_{\geq 0}} \HH^m(A, \Omega_{nc}^p(A)).$$
%Let $A$ be a radical square zero algebra with a Wedderbrun-Malcev decomposition $A=E\oplus r$.
%Let $\Omega_{nc}^p(A)$ be the $p$-th kernel in the projective resolution $\calR(A)$.Given $f\in \Hom_{E-E}(r^{\otimes_E m}, \Omega_{nc}^p(A))$ and $g\in\Hom_{E-E}(r^{\otimes_E n}, \Omega_{nc}^{q}(A)),$ the Gerstenhaber bracket $[\cdot, \cdot]$ can be written as
%\begin{equation}\label{equa-ger}
%\begin{split}
%  [f, g](a_{1, m+n-1})=&\sum_{i=1}^m (-1)^{(i-1)(q-n-1)}(f\otimes \id)(a_{1, i-1}\otimes
%  g(a_{i, i+n-1})\otimes a_{i+n, m+n-1})+\\
%  &\sum_{i=1}^q(-1)^{i(p-m-1)}(\id^{\otimes i}\otimes f)(g(a_{1, n})\otimes a_{n+1, m+n-1})-
%  (-1)^{(m-p-1)(n-q-1)}\\
%  &\sum_{i=1}^{n}(-1)^{(i-1)(p-m-1)}(g\otimes \id)(a_{1, i-1}\otimes f(a_{i, i+m-1})\otimes a_{i+m, m+n-1})-\\
%  &(-1)^{(m-p-1)(n-q-1)}\sum_{i=1}^p(-1)^{i(q-n-1)}(\id^{\otimes i}\otimes g)(f(a_{1,m})\otimes a_{m+1, m+n-1}),
%\end{split}
%\end{equation}
%for any $a_{1, m+n-1}\in r^{\otimes_E m+n-1}$.
In particular, for $m=n=1$, we have
\begin{equation}\label{equa-ger2}
  \begin{split}
    [f, g](a)=& \sum_{i=1}^q(-1)^{ip}(\id^{\otimes i}\otimes f)(g(a))-(-1)^{pq}\sum_{i=0}^p(-1)^{iq}
    (\id^{\otimes i}\otimes g)(f(a)).
  \end{split}
\end{equation}

Let us introduce the notation $\diamond$ (cf. \cite{San, San2}):
Given two paths $\alpha\in Q_m$  and $\beta\in Q_n$. Suppose that
$\alpha=a_1a_2\cdots a_m,$ and $\beta=b_1b_2\cdots b_n.$
where $a_i, b_j\in Q_1$. Let $i=1, \cdots, m. $  If $a_i//\beta,$ we denote by $\alpha\diamond_i\beta$ the path  in $Q_{m+n-1}$
 obtained by replacing the arrow $a_i$ by the path $\beta$. Namely, we define
 \begin{equation*}
 \alpha\diamond_i\beta:=
 \begin{cases}
 a_1\cdots a_{i-1}b_1\cdots b_na_{i+1}\cdots a_{m}, & \mbox{if $a_i// \beta$},\\
 0 & \mbox{otherwise.}
 \end{cases}
 \end{equation*}
\begin{lemma}\label{lemma-bracket-quiver}
Let $k$ be a field. Let $Q$ be a finite and connected quiver. Then via the linear isomorphism in Lemma \ref{lemma-new},
we have the following
\begin{enumerate}
\item Let $(x, y)\in k(Q_1//Q_1)\subset\Hom_{E-E}(r, A)$ and $$(\gamma_m, \beta_{p+1})\in k(Q_{m}//Q_{p+1})\subset\Hom_{E-E}(r^{\otimes_E m}, \Omega_{nc}^{p}(A)).$$
Then $$[(x, y), (\gamma_m, \beta_{p+1})]=\sum_{i=1}^{p+1}\delta_{b_i, x}(\gamma_m, \beta_{p+1}\diamond_iy)-\sum_{i=1}^m\delta_{a_i, y}(\gamma_{m}\diamond_i x, \beta_{p+1}),$$
where $a_i, b_i$ are the $i$-th arrow in $\gamma_m$ and $\beta_{p+1}$,
respectively.
\item Let $(x, \gamma_{p+1})\in k(Q_1//Q_{p+1})\subset\Hom_{E-E}(r, \Omega_{nc}^p(A))$ and $$(y, \beta_{q+1})\in k(Q_1//Q_{q+1})\subset\Hom_{E-E}(r, \Omega_{nc}^q(A)).$$
Then we have 
\begin{equation*}
\begin{split}
[(x, \gamma_{p+1}), (y, \beta_{q+1})]=&\sum_{i=1}^{q+1}(-1)^{(i-1)p}\delta_{x, b_i}(y, \beta_{q+1}\diamond_i \gamma_{p+1})-\\
&(-1)^{pq}\sum_{i=1}^{p+1}(-1)^{(i-1)q}\delta_{y, a_i} (x, \gamma_{p+1}\diamond_i \beta_{q+1}),
\end{split}
\end{equation*}
where $a_i, b_i$ are the $i$-th arrow in $\gamma_{p+1}$ and $\beta_{q+1}$,
respectively.
\end{enumerate}
\end{lemma}
\pf This is a direct consequence of Formula (\ref{equa-ger2}).
\begin{rem}
For the general case, the formula of the Gerstenhaber bracket $[\cdot, \cdot]$ is quite complicated. We will use prop theory (e.g., \cite{LoVa, Mar}) to describe it in a future paper. %n the following, we will use Lemma \ref{lemma-bracket-quiver} to compute the Gerstenhaber bracket on $\HH_{\sg}^*(A, A)$ for $c$-crowns. 
\end{rem}

%In the following sections, we will focus on some concrete examples of quivers.
\section{$c$-crown}
In this section, we will use Lemma \ref{lemma-bracket-quiver} to  compute the Gerstenhaber algebra and BV algebra structures on $\HH_{\sg}^*(A, A)$  of the radical square zero algebras $A$ associated with $c$-crown quivers (cf. Definition \ref{defn-crown}).  More precisely, we prove that the Gerstenhaber algebra  $(\HH^*(A, A), \cup, [\cdot, \cdot])$ for the $c$-crown quivers with $c\in 2\mathbb Z$,  is isomorphic to the semidirect product of the Witt algebra $\mathcal W$ and the Laurent polynomial ring $\mathcal M=k[t, t^{-1}]$, after some grading shifts (cf. Theorem \ref{equ-bv00}). Throughout this section, we assume that the base field $k$ is not of characteristic two. 

\subsection{The case: $c=1$}
Let us first consider the $1$-crown (i.e. the one loop quiver).
$$
\begin{tikzcd}
Q:=\bullet \arrow[loop right]{r}{a}
\end{tikzcd}$$
Its radical square zero algebra is $A=k[a]/(a^2)$, the algebra of
dual numbers. Since $A$ is a commutative symmetric algebra, from   \cite[Corollary
6.4.1]{Bu}, we have that
\begin{equation*}
  \HH_{\sg}^m(A, A)=
  \begin{cases}
    \HH^m(A, A) & \mbox{for $m>0$,}\\
    \Tor_{-m-1}^{A^e}(A, A) & \mbox{for $m<-1$.}
  \end{cases}
\end{equation*}
\begin{prop}[Proposition 3.4 \cite{Cib3}]\label{prop-cib3}
  Let $Q$ be the one loop quiver and $A$ be its radical square zero algebra over a field $k$. Assume that $k$ is not of characteristic two. Then for every $n>0$, we have $\dim \HH^n(A, A)=1.$
\end{prop}
\begin{rem}\label{rem-cib3}
Since $A$ is a symmetric algebra, we have $\HH^m(A, A)\cong \Tor_m^{A^e}(A, A)^*$. Thus $\dim \Tor_m^{A^e}(A, A)=1$ for $m>0$. 
\end{rem}
\begin{lemma}
For any $n\in \Z$, we have $\dim \HH_{\sg}^n(A, A)=1.$
\end{lemma}
\pf From Proposition \ref{prop-cib3} and Remark \ref{rem-cib3},  it is sufficient to show that
$$\dim \HH_{\sg}^0(A, A)=\dim \HH_{\sg}^{-1}(A, A)=1.$$
Recall that from \cite[Corollary
6.4.1]{Bu}, we have an exact sequence,
\begin{equation*}
\xymatrix{
0 \ar[r] & \HH_{\sg}^{-1}(A, A)\ar[r] & \HH_0(A, A)\ar[d]^-{\cong} \ar[r] &  \HH^0(A, A)\ar[d]^-{\cong} \ar[r]  & \HH_{\sg}^0(A, A)\ar[r]  & 0\\
&   & k[a]/(a^2)\ar[r]^-{\mu_a} & k[a]/(a^2)
}
\end{equation*}
where $\mu_a$ is the map multiplied by $a$. This follows that  $$\dim \HH_{\sg}^{-1}(A, A)=\dim \HH_{\sg}^0(A, A)=1.$$ \epf
\begin{rem}
The Lie algebra structure on $\HH^*(A, A)$ has been investigated in \cite{San}. Next we will describe the BV algebra structure on $\HH_{\sg}^*(A, A)$.
\end{rem}
Since $A$ is symmetric,  we have for $m\in\Z_{\geq 0}$, $$\HH_{\sg}^{-m}(A, A)\cong \HH^1(A, \Omega_{nc}^{m+1}(A)).$$ From Remark \ref{rem-coho}, it follows that $$\HH^1(A, \Omega_{nc}^{m+1}(A))\cong  \frac{k(Q_1//Q_{m+2})}{\Imm(D_{0, m+1})}\oplus\Ker(D_{1, m+1}).$$ Recall that $D_{i, m+1}: k(Q_i//Q_{m+1})\rightarrow k(Q_{i+1}//Q_{m+2})$ is defined as follows,
$$D_{i, m+1}(\gamma_i, \beta_{m+1}):=\sum_{a\in Q_1} (a\gamma_i, a\beta_{m+1})+(-1)^{i+m}\sum_{a\in Q_1}(\gamma_i
a, \beta_{m+1}a).$$
Note that if $m$ is odd, then $D_{0, m+1}=0$ and $D_{1, m+1}$ is a bijection. Similarly if $m$ is even, then $D_{0, m+1}$ is a bijection and $D_{1, m+1}=0.$ So we have
\begin{equation}\label{equ-isom3}
  \HH^1(A, \Omega_{nc}^{m+1}(A))=
  \begin{cases}
    k(Q_1//Q_{m+2}) & \mbox{if $m$ is odd},\\
    k(Q_1//Q_{m+1}) & \mbox{if $m$ is even}.
  \end{cases}
\end{equation}
%Recall that, for a radical square zero algebra $A$, the Gerstenhaber bracket $[\cdot,\cdot]$ on $\HH^*_{\sg}(A, A)$  can be written as follows:
%Let $f\in \HH^1(A, \Omega_{nc}^p(A))$ and $g\in \HH^1(A, \Omega_{nc}^q(A))$,
%then for  any $a\in \rad(A)$,
%\begin{equation*}
%  [f, g](a):=\sum_{i=0}^q(-1)^{ip}(\id_i\otimes f)(g(a))-(-1)^{pq}\sum_{i=1}^p
%  (-1)^{iq}(\id_i\otimes g)(f(a)).
%\end{equation*}
%From this formula, we can describe the Gerstenhaber bracket $[\cdot,\cdot]$ under the isomorphisms (\ref{equ-isom3}).
\begin{prop}\label{prop-Gerb}
  Let $Q$ be the one loop quiver
 and $A$ be its radical square zero algebra over
  a field $k$. Assume that $k$ is not of characteristic two.  Then for $m, n\in\Z_{\geq 1}$, we have the following cases:
  \begin{enumerate}
    \item If both $m$ and $n$ are odd, then we have the following commutative diagram,
    \begin{equation*}
    \xymatrix{
      \HH^1(A,\Omega_{nc}^{m+1}(A))\times \HH^1(A, \Omega_{nc}^{n+1}(A))\ar[d]^{\cong}\ar[r]^-{[\cdot, \cdot]} & \HH^1(A, \Omega_{nc}^{m+n+2}(A))\ar[d]^-{\cong}\\
      k(Q_1//Q_{m+2}) \times k(Q_1//Q_{n+2} \ar[r]^-{\{\cdot, \cdot\}}) & k(Q_1//Q_{m+n+3})
      }
    \end{equation*}
    where the bracket $\{\cdot,\cdot\}$ is defined as follows,
    $$\{(a, a^{m+2}), (a, a^{n+2})\}=(n-m)(a, a^{m+n+3}).$$
    \item If both $m$ and $n$ are even, then we have the following commutative diagram,
    \begin{equation*}
    \xymatrix{
      \HH^1(A,\Omega_{nc}^{m+1}(A))\times \HH^1(A, \Omega_{nc}^{n+1}(A))\ar[d]^{\cong}\ar[r]^-{[\cdot, \cdot]} & \HH^1(A, \Omega_{nc}^{m+n+2}(A))\ar[d]^-{\cong}\\
      k(Q_1//Q_{m+1}) \times k(Q_1//Q_{n+1} \ar[r]^-{\{\cdot, \cdot\}}) & k(Q_1//Q_{m+n+3})
      }
    \end{equation*}
    where $\{(a, a^{m+1}), (a, a^{n+1})\}=0.$
    \item If $m$ is even and $n$ is odd, then the following diagram commutes,
    \begin{equation*}
    \xymatrix{
      \HH^1(A,\Omega_{nc}^{m+1}(A))\times \HH^1(A, \Omega_{nc}^{n+1}(A))\ar[d]^{\cong}\ar[r]^-{[\cdot, \cdot]} & \HH^1(A, \Omega_{nc}^{m+n+2}(A))\ar[d]^-{\cong}\\
      k(Q_1//Q_{m+1}) \times k(Q_1//Q_{n+2} \ar[r]^-{\{\cdot, \cdot\}}) & k(Q_1//Q_{m+n+2})
      }
    \end{equation*}
    where $\{(a, a^{m+1}), (a, a^{n+2})\}:=-m(a, a^{m+n+2}).$
  \end{enumerate}
\end{prop}
\pf The assertion (1) comes from Lemma \ref{lemma-bracket-quiver}. Now let us prove the assertion (2). Recall that $(a, a^{m+1})$ represents the element in $\Hom_{E-E}(r, \Omega_{nc}^{m+1}(A))$, which sends $a$ to $ea^{m+1}+(-1)^{m+1}a^{m+1}e\in \Omega_{nc}^{m+1}(A).$ From Formula (\ref{equa-ger2}) it follows that
$$[(a, a^{m+1}), (a, a^{n+1})]\in k(Q_1//Q_{m+n+2})\subset \Hom_{E-E}(r, \Omega_{nc}^{m+n+2}(A)).$$ So from Formula (\ref{equ-isom3}), we have $[(a, a^{m+1}), (a, a^{n+1})]=0$ in $\HH^1(A, \Omega_{nc}^{m+n+2}(A))$.
It remains to verify the assertion (3).
From Formula (\ref{equa-ger2}) again, we have
\begin{equation*}
\begin{split}
[(a, a^{m+1}), (a, a^{n+2})](a)=&\sum_{i=0}^{n+1}(-1)^{i(m+1)}(\id^{\otimes i}\otimes (a, ea^{m+1}+(-1)^{m+1}a^{m+1}e))(a^{n+2})-\\
&\sum_{i=0}^{m+1}(-1)^{(i+m+1)(n+1)} (\id^{\otimes i}\otimes (a, a^{n+2}))(ea^{m+1}+(-1)^{m+1}a^{m+1}e)\\
%=&\sum_{i=0}^{n+1}(-1)^{i}(\id^{\otimes i}\otimes (ea^{m+1}-a^{m+1}e))(a^{n+2})-\\
%&\sum_{i=0}^{m+1} (\id^{\otimes i}\otimes (a, a^{n+2}))(ea^{m+1}-a^{m+1}e)\\
%&=ea^{m+n+2}-a^{m+n+2}e-(m+1)(ea^{m+n+2}-a^{m+n+2}e)\\
&=-m(ea^{m+n+2}-a^{m+n+2}e),
\end{split}
\end{equation*}
thus $[(a, a^{m+1}), (a, a^{n+2})]=-m(a, a^{m+n+2})\in \HH^1(A, \Omega_{nc}^{m+n+2}(A)).$ This completes the proof.
\epf

From  \cite[Theorem 6.17]{Wang} it follows that
$(\HH_{\sg}^{*}(A, A), \cup, [\cdot, \cdot], \Delta)$ is a BV algebra. Now let us describe the BV algebra structure in this concrete example.
First, we can write down the formula for the Connes B-operator.
\begin{lemma}\label{lemma-bv-1}
  Let $Q$ be the one loop quiver and $A$ be its radical square zero algebra.
  Then for $m\in\Z_{\geq 0}$, we have
  \begin{enumerate}
    \item If $m$ is even, then we have the following commutative diagram,
    \begin{equation*}
      \xymatrix{
      \HH^1(A, \Omega_{nc}^{m+1}(A))\ar[d]^{\cong} \ar[r]^-{B}& \HH^1(A, \Omega_{nc}^{m+2}(A))\ar[d]^{\cong}\\
      k(Q_1//Q_{m+1}) \ar[r]^{\Delta=0} & k(Q_1// Q_{m+3})
      }
    \end{equation*}
    \item If $m$ is odd, then the following diagram commutes,
    \begin{equation*}
      \xymatrix{
      \HH^1(A, \Omega_{nc}^{m+1}(A))\ar[d]^{\cong} \ar[r]^-{B}& \HH^1(A, \Omega_{nc}^{m+2}(A))\ar[d]^{\cong}\\
      k(Q_1//Q_{m+2}) \ar[r]^{\Delta} & k(Q_1// Q_{m+2}),
      }
    \end{equation*}
    where $\Delta$ is defined as $\Delta( (a, a^{m+2}))=m(a, a^{m+2}).$
  \end{enumerate}
\end{lemma}
\pf The proof is completely analogous to the one of Proposition \ref{prop-Gerb}.
\epf

Similarly, we can also write down the formula for the cup product $\cup$.
\begin{lemma}\label{lemma-cup-1}
   Let $Q$ be the one loop quiver and $A$ be its radical square zero $k$-algebra, where $k$
   is not of characteristic two.
   Then we have the following cases for $m, n\in \Z_{>0}$,
   \begin{enumerate}
     \item If both $m$ and $n$ are odd, then we have the following commutative diagram,
     \begin{equation*}
       \xymatrix{
      \HH^1(A,\Omega_{nc}^{m+1}(A))\times \HH^1(A, \Omega_{nc}^{n+1}(A))\ar[d]^{\cong}\ar[r]^-{\cup} & \HH^1(A, \Omega_{nc}^{m+n+1})\ar[d]^-{\cong}\\
      k(Q_1//Q_{m+2}) \times k(Q_1//Q_{n+2}) \ar[r]^-{0} & k(Q_1//Q_{m+n+1})
      }
     \end{equation*}
     where we used the connecting isomorphism $$\theta_{1, m+n+1}:\HH^1(A, \Omega_{nc}^{m+n+1}(A))
     \rightarrow \HH^2(A, \Omega_{nc}^{m+n+2}(A)).$$
     \item If both $m$ and $n$ are even, then
     \begin{equation*}
       \xymatrix{
      \HH^1(A,\Omega_{nc}^{m+1}(A))\times \HH^1(A, \Omega_{nc}^{n+1}(A))\ar[d]^{\cong}\ar[r]^-{\cup} & \HH^1(A, \Omega_{nc}^{m+n+1})\ar[d]^-{\cong}\\
      k(Q_1//Q_{m+1}) \times k(Q_1//Q_{n+1}) \ar[r]^-{\cup'} & k(Q_1//Q_{m+n+1})
      }
     \end{equation*}
     where $(a, a^{m+1})\cup'(a, a^{n+1})=(a, a^{m+n+1}).$
     \item If $m$ is even and $n$ is odd, then
     \begin{equation*}
       \xymatrix{
      \HH^1(A,\Omega_{nc}^{m+1}(A))\times \HH^1(A, \Omega_{nc}^{n+1}(A))\ar[d]^{\cong}\ar[r]^-{\cup} & \HH^1(A, \Omega_{nc}^{m+n+1})\ar[d]^-{\cong}\\
      k(Q_1//Q_{m+1}) \times k(Q_1//Q_{n+2}) \ar[r]^-{\cup'} & k(Q_1//Q_{m+n+2})
      }
     \end{equation*}
     where $(a, a^{m+1})\cup'(a, a^{n+2})=-(a, a^{m+n+2}).$
   \end{enumerate}
\end{lemma}
\pf The proof is completely analogous to the one of Proposition \ref{prop-Gerb}.
\epf

\begin{rem}\label{rem-witt}
 % Let us recall the Witt algebra $\calW$ (cf. e.g. \cite{Zim}). As a vector space, $\calW=\bigoplus_{n\in\Z}k\langle L_n\rangle,$ where $k\langle L_n\rangle$ is a one-dimensional $k$-vector space with a basis $L_n.$ The Lie bracket is defined  $$[L_m, L_n]=(m-n)L_{m+n}.$$ Clearly, the even part of $\calW^{}$ $$\calW^{even}:=\bigoplus_{m\in \Z} k\langle L_{2m}\rangle$$ is a Lie subalgebra of $\calW$. Let us construct a natural representation $\calM$ of $\calW$. As a vector space, $$\calM=\bigoplus_{n\in\Z} k\langle M_n\rangle,$$ where $k\langle M_n\rangle $ is a one-dimensional $k$-vector space with a basis $M_n$. Define the action of $L_m$ on $M_n$ as follows: $$[L_m, M_n]:=-nM_{m+n}.$$ Then one can check that it induces a representation on $\calM$ of $\calW$. Denote $$\calM^{even}:=\bigoplus_{n\in\Z}k\langle M_{2n+1} \rangle.$$ Trivialy $\calM^{even}$ is a representation of $\calW^{even}$.
Let  $\calW$  be  the Lie algebra of vector fields with Laurent polynomial coefficients, i.e. those of the form $f(t)\frac{d}{dt},$ with $f(t)\in k[t, t^{-1}]$. The Lie bracket is given by $$[f(t)\frac{d}{dt}, g(t)\frac{d}{dt}]=f(t)\frac{g(t)}{dt} \frac{d}{dt}-g(t)\frac{f(t)}{dt} \frac{d}{dt}.$$ %Then $$L_n:=-t^{n+1}\frac{d}{dt}$$ for $n\in \Z$ is a basis of $\calW$. It is straightforward to verify that $$[L_m, L_n]=(m-n)L_{m+n}.$$
Denote by $\calM$ the Laurent polynomial ring $k[t, t^{-1}]$. %Clearly it has a basis $$M_n:=t^n$$ for $n\in \Z$ and 
Clearly, $\calW$ acts on $\calM$ by derivations, namely,  for any $g(t)\in k[t, t^{-1}]$, we define $[f(t)\frac{d}{dt}, g(t)]:= f(t) \frac {dg(t)}{dt}.$
  %Hence we have, $$[L_m, M_n]=-nt^{m+n}=-nM_{m+n}.$$
  It is straightforward to verify that this action defines a Lie module structure on $\calM$ over the Lie algebra $\calW$. Note that we also have an action of the commutative algebra $\calM$ on $\calW$:
  For $f(t)\in \calM$ and $g(t)\frac{d}{dt}\in \calW$, $$f(t)\cdot g(t)\frac{d}{dt}:=f(t)g(t)\frac{d}{dt}\in \calW.$$
  %In particular, we have $$M_n\cdot L_m=L_{m+n}.$$
  %Clearly, by this action $\calW$ is a module over the commutative algebra $\calM$.
  
Let us construct a BV algebra $(\calM\times \calW, \cup,
  [\cdot, \cdot], \Delta)$ as follows.
  The grading is given by
  \begin{equation*}
    (\calM\times \calW)_n:=
    \begin{cases}
      k\langle M_m \rangle & \mbox{if $n=2m$;}\\
      k\langle L_m\rangle  & \mbox{if $n=2m+1$,}
    \end{cases}
  \end{equation*}
  where $M_m=t^m$ and $L_m=-t^{m+1}\frac{d}{dt}$.
  As a graded commutative algebra, $(\calM\times \calW, \cup)=(\calM \ltimes \calW, \cdot)$  and as a graded Lie algebra, $(\calM\times \calW, [\cdot, \cdot])= (\calM\rtimes \calW, [\cdot, \cdot]).$
%We define a graded commutative algebra structure on it:
%\begin{equation*}
%\begin{split}
%M_{2m}\cup M_{2n}:&=M_{2m+2n},\\
%M_{2m}\cup L_{2n}:&=-L_{2m+2n},\\
%L_{2m}\cup L_{2n}:&=0.
%\end{split}
%\end{equation*}
The BV operator is defined as
\begin{equation*}
\begin{split}
\Delta_{2m}(M_{m}):&=0,\\
\Delta_{2m+1}(L_{m}):&=-mM_{m}.
\end{split}
\end{equation*}
Then one can check that $(\calM\times \calW, \cup, [\cdot,\cdot], \Delta)$
is a BV algebra.
Similarly, we can also construct a BV algebra $(\calM^{even}\times \calW^{even},
\cup, [\cdot, \cdot], \Delta)$.
The grading is
\begin{equation*}
(\calM^{even}\times \calW^{even})_n=
\begin{cases}
k\langle L_{n-1}\rangle  & \mbox{if $n$ is odd},\\
k\langle M_{n}\rangle & \mbox{if $n$ is even}.
\end{cases}
\end{equation*}
As a graded commutative algebra, $(\calM^{even}\times \calW^{even}, \cup)\cong(\calM^{even}\ltimes \calW^{even}, \cdot)$ and as a graded Lie algebra, $(\calM^{even}\times \calW^{even}, [\cdot, \cdot])\cong
\calM^{even}\rtimes  (\calW^{even}, [\cdot, \cdot]).$
The BV operator is defined by
\begin{equation*}
\begin{split}
\Delta_{2m}(M_{2m}):&=0,\\
\Delta_{2m+1}(L_{2m}):&=-2mM_{2m}.
\end{split}
\end{equation*}

 %and it is isomorphic to the BV algebra $(\HH_{\sg}^*(A, A), \cup, \Delta)$

\end{rem}

%Now let us describe the Gerstenhaber algebra structure on $\HH^*_{\sg}(A, A)$. Denote, for $m\in \Z$,
%\begin{equation*}
%  \LL_{2m}:=\HH^{2m+1}_{\sg}(A, A)
%\end{equation*}
%and
%\begin{equation*}
% \MM_{2m}:= \HH_{\sg}^{2m}(A, A).
%  \end{equation*}
%Then we have the following proposition.
\begin{prop}\label{prop-de-ger}
   Let $Q$ be the one loop quiver
 and $A$ be its radical square zero algebra over
  a field $k$, where $k$ is not of characteristic two. Then %we have a Lie algebra isomorphism $$(\HH_{\sg}^{odd}(A, A)[1], [\cdot, \cdot])\cong(\calW^{even}, [\cdot, \cdot]).$$ and an isomorphism of commutative algebras $$(\HH_{\sg}^{even}(A, A), \cup)\cong (\calM^{even}, \cdot).$$
$(\HH_{\sg}^*(A, A), \cup, [\cdot, \cdot], \Delta)$ is isomorphic to the BV algebra $(\calM^{even}\times \calW^{even}, \cdot, [\cdot,\cdot], \Delta).$
\end{prop}
\pf It is a direct consequence of Proposition \ref{prop-Gerb}, Lemmas \ref{lemma-bv-1} and
\ref{lemma-cup-1}.
\epf

%Hence, as a conclusion of  the two lemmas above, we obtain the following proposition.
%\begin{prop}
%  Let $Q$ be the one loop quiver and $A$ be its radical square zero algebra.
%  Then $\HH_{\sg}^{*}(A, A)$ is a BV algebra equipped with the Connes B-operator.
%  Hence, it induces a BV algebra structure on the semi-direct product Lie algebra $\calW^{even}\ltimes \calM^{even}$.
%\end{prop}
%Now let us describe this BV algebra structure on $\calW^{even}\ltimes \calM^{even}$. As a graded space,

\subsection{The case: $c\geq 2$}
In this section, we fix $Q$ to be a $c$-crown with $c\geq 2$.
Denote by $A$ its radical square zero $k$-algebra.  We assume that  $k$ is not of characteristic two. Let us first recall a result in \cite{Cib3}.
\begin{prop}[Proposition 3.3 \cite{Cib3}]\label{prop3.3}
Let $Q$ be a $c$-crown with $c\geq 2$. Let $n$ be an even multiple of $c$. Then $\dim_k \HH^n(A, A)=\dim_k\HH^{n+1}(A, A)=1.$
The cohomology vanishes in all other degrees. 
\end{prop}
Next let us consider the Tate-Hochschild cohomology $\HH^*_{\sg}(A, A)$. Note that $A$ is a self-injective algebra (but not a symmetric algebra). From \cite[Corollary 6.4.1]{Bu}, it follows that for $m\geq 1,$ $\HH^m_{\sg}(A, A)\cong \HH^m(A, A)$ and that for $m>1$,  $\HH_{\sg}^{-m}(A, A)\cong \HH^1(A, \Omega_{nc}^{m+1}(A)).$ By Remark \ref{rem-coho},  $ \HH^1(A, \Omega_{nc}^{m+1}(A))\cong \frac{k(Q_1//Q_{m+2})}{\Imm(D_{0, m+1})}\oplus \Ker(D_{1, m+1}).$ Recall that
\begin{equation*}
  \begin{split}
    D_{0, m+1}(e, \gamma_{m+1})&=\sum_{a\in Q_1}(a, a\gamma_{m+1})
    +(-1)^m\sum_{a\in Q_1}(a, \gamma_{m+1}a),\\
    D_{1, m+1}(x, \gamma_{m+1})&=\sum_{a\in Q_1}(ax, a\gamma_{m+1})+
    (-1)^{m+1}\sum_{a\in Q_1}(xa, \gamma_{m+1}a).
  \end{split}
\end{equation*}
Anaglous to the proof of Proposition \ref{prop3.3} above, we have the following result.
\begin{prop}
  Let $Q$ be a $c$-crown with $c\geq 2$. If $m$ is an even multiple of $c$,
  then $$\dim_k \HH^{m}_{\sg}(A, A)=\dim_k\HH^{m+1}_{\sg}(A, A)=1.$$
 The Tate-Hochschild cohomology vanishes in all other degrees.
\end{prop}

\begin{prop}
  Let $Q$ be a $c$-crown with $c\in 2\Z_{>0}$ and $A$ be its radical square zero algebra over a field $k$, where $k$ is not of characteristic two. Denote by $\gamma$ the oriented cycle (i.e. $\gamma=a_1a_2\cdots a_c$).
  Then we have the following cases:
  \begin{enumerate}
    \item Let $(a, a\gamma^p)\in \HH^1(A, \Omega_{nc}^{cp}(A))$ and $(a, a\gamma^q)\in \HH^1(A, \Omega_{nc}^{cq}(A))$ be nontrivial elements respectively, where $a\in Q_1$.
    Then $$[(a, a\gamma^p), (a, a\gamma^q)]=(q-p) (a, a\gamma^{p+q}).$$
    \item Let $x:=\sum_{a\in Q_1}(a, a\gamma^p)\in \HH^1(A, \Omega_{nc}^{cp+1}(A))$and $$y:=\sum_{a\in Q_1}(a, a\gamma^q)\in \HH^1(A, \Omega_{nc}^{cq+1})$$ be nontrivial elements, respectively,
    then $[x, y]=0.$
    \item Let $x:=\sum_{a\in Q_1}(a, a\gamma^p)\in \HH^1(A, \Omega_{nc}^{cp+1}(A))$ and $$(b, b\gamma^q)\in\HH^1(A, \Omega_{nc}^{cq}(A))$$ be nontrivial elements, respectively,
        then $$[x, (b, b\gamma^q)]=-p\sum_{a\in Q_1}(a, a\gamma^{p+q}).$$
  \end{enumerate}
\end{prop}
%Denote, for $p\in\Z_{\geq 0}$,
%$$L_p:=\HH^1(A, \Omega_{nc}^{cp}(A))$$ and $$M_p:=\HH^1(A, \Omega_{nc}^{cp+1}(A)).$$
%Then we have that  $$L_{\bullet}:=\bigoplus_{p=0}^{\infty} L_p$$ is an opposite Witt algebra and
%$$M_{\bullet}:=\bigoplus_{p=0}^{\infty} M_p$$ is a $L_{\bullet}$-module.
%Denote, for $p\in \Z_{\geq 0}$,
%$$L^p:=\HH^{cp+1}(A, A)$$
%and
%$$M^p:=\HH^{cp}(A, A).$$
%Similarly, form \cite{San}, we obtain that
%$$L^{\bullet}:=\bigoplus_{p=0}^{\infty} L^p$$
%is a Witt algebra and $$M^{\bullet}:=\bigoplus_{p=0}^{\infty} M^p$$
%is a $L^{\bullet}$-module. Now we are interested in the Gerstenhaber Lie bracket
%between $(L_{\bullet}, M_{\bullet})$ and $(L^{\bullet}, M^{\bulle}t})$.
%Then we have the following proposition.
%\begin{prop}
%  $[M^{\bullet}, M_{\bullet}]=0.$
%\end{prop}

For $p\in \Z$, we denote $ \mathbb{L}_p:=\HH_{\sg}^{cp+1}(A, A)$ and $\MM_p:=\HH_{\sg}^{cp}(A, A).$

%Then we have the following proposition.
\begin{thm}\label{equ-bv00}
  Let $Q$ be a $c$-crown with $c\in2\Z_{>0}$ and $A$ be its radical
  square zero $k$-algebra. Assume that $k$ is not of characteristic two. Then $(\LL:=\bigoplus_{i\in\Z} \LL_i, [\cdot, \cdot])$ is isomorphic to the Witt algebra $\calW$ and $(\MM:=\bigoplus_{i\in \Z}\MM_i, \cup)$ is isomorphic to the graded commutative algebra  $\calM$. %As a consequence,  $(\HH_{\sg}^*(A, A), \cup, [\cdot, \cdot])$ is isomorphic to $(\calM^{even}\times \calW^{even},\cup, [\cdot, \cdot], \Delta)$. 
  Moreover, $(\HH^*_{\sg}(A, A), \cup, [\cdot, \cdot])$
 is isomorphic to $(\calM\times \calW, \cup, [\cdot, \cdot])$ as Gerstenhaber algebras.
\end{thm}
\pf The proof is completely analogous to the one of Proposition \ref{prop-de-ger}.
\epf

\section{$r$-multiple loops quiver}\label{two loops quiver}
In this section, we study   the Gerstenhaber algebra structure (cf. Theorem \ref{prop-prop-two}) on the Tate-Hochschild cohomology of the radical square zero algebra associated to the  $r$-multiple loops quiver $Q_r$ ($r\geq 2$).
  $$
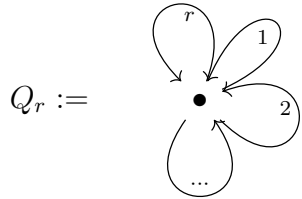
\begin{tikzcd}
Q_r:= & \bullet \ar[%
    ,loop % tells tikz-cd to do a loop
    ,out=25 % start at angle 123?
    ,in=67 % stop at angle 57?
    ,distance=4em % biggest distance of arrow to node. Yarou can use pt or cm as well.
    ]{}{1}
    \ar[
    ,loop
    ,out=70
    ,in=135
    ,distance=4em
    ]{}{r}
     \ar[
    , loop
    ,in=20
    ,out=315
    ,distance=4em]{}{2}
    \ar[
    ,loop
    ,in=305
    ,out=235
    ,distance=4em]{}{\cdots}
\end{tikzcd}
$$
%$$
%\begin{tikzcd}
%Q:=& \bullet  \arrow[loop right]{r}{b}  \arrow[loop left]{l}{a}
%\end{tikzcd}$$
We denote the radical square zero algebra by $A\cong k[x_1, \cdots, x_r]/(\{x_ix_j \}_{1\leq i, j\leq r})$ over a field $k$ of characteristic zero.
\begin{prop}[Proposition 4.4. \cite{San2}]\label{prop-san}
Assume that $Q$ is the two loops quiver with the loops $1$ and $2$. Let $A$ be its radical square zero algebra over a field $k$ of characteristic zero. Then $\HH^1(A, A)\cong k(Q_1//Q_1)$ and the elements
in $k(Q_1// Q_1)$
\begin{equation*}
\begin{split}
H:=&(1, 1)-(2, 2),\\
E:=&(2,1),\\
F:=&(1, 2)
\end{split}
\end{equation*}
generate a copy of the Lie algebra $\mathfrak{sl}_2(k)$ in $\HH^1(A, A)$. Moreover, the Lie algebra $\HH^1(A, A)$ is isomorphic to $\Sl_2\times k$, where $I:=(1, 1)+(2, 2)$ is a non-zero element such that $[I, \HH^1(A, A)]=0$.
\end{prop}
\begin{rem}
%If we denote the path $a$ by 1 and the path $b$ by 2. We also denote $(x, y):=E_{y, x},$ where $\{x, y\}\subset\{1, 2\}$. Then we have
%\begin{equation*}
%\begin{split}
%H=&E_{1, 1}-E_{2, 2}\\
%F=&E_{2, 1},
%\end{split}
%\end{equation*}
%where $E_{i, j}$ is the basis of $\Gl_2(k)$, which is defined as follows,
%$$E_{i, j}(e_k)=\delta_{i, k} e_j.$
The author in \cite{San2}  describes the module structure $\HH^m(A, A)$ over $\HH^1(A, A)$ for any $m\in \Z_{>0}$. Next we will completely determine the Gerstenhaber algebra structure on $\HH^*_{\sg}(A, A)$.
\end{rem}

Let us first work on a more general setting. Let $V$ be a $k$-linear vector space of dimension $r$.  For $m, p\in \Z_{> 0}$, denote $$T^{m, p}(V):=T_0^{m, p}(V)\oplus T_1^{m, p}(V),$$ where $T_0^{m, p}(V):=\Hom_k(V^{\otimes m},V^{\otimes p})$
and $T_1^{m, p}(V):=\Hom_{k}(V^{\otimes m}, V^{\otimes p+1}).$
We denote
$$T^{*, *}(V):=\bigoplus_{m>0,  p\geq 0} T^{m, p}(V).$$
%where we use the notation $V^{\otimes 0}:=k$.
%Note that $T^{>0}(\End(V))$ is a subspace of $T^{*, *}(V)$ and
%the bracket $\{\cdot,\cdot\}$ on $T^{>0}(\End(V))$ can extend to
%$T^{*, *}(V)$.
We define a bracket $\{\cdot, \cdot\}$ on $T^{*, *}(V)$ as follows:
Let $f\in T_1^{m, p}(V)$ and $g\in T_1^{n, q}(V)$, we define
$\{f, g\}\in T_1^{m+n-1, p+q}(V)$ as 
\begin{equation}\label{defn-bracket-T}
  \begin{split}
    \{f, g\}:=&\sum_{i=1}^m(-1)^{(i-1)(q-n-1)}(f\otimes \id^{\otimes q})(\id^{\otimes i-1}\otimes g\otimes \id^{\otimes m-i})+\\
    &\sum_{i=1}^{q-1}(-1)^{i(p-m-1)}(\id^{\otimes i}\otimes f\otimes \id^{\otimes q-i})(g\otimes \id^{\otimes m-1})-\\
    &\sum_{i=1}^n(-1)^{(n-q+i)(p-m-1)}(g\otimes \id^{\otimes p})(\id^{\otimes i-1}
    \otimes f\otimes \id^{\otimes n-i})-\\
    &\sum_{i=1}^{p-1}(-1)^{(m-p+i-1)(q-n-1)}(\id^{\otimes i}
    \otimes g\otimes\id^{\otimes p-i})(f\otimes \id^{\otimes n-1}).
  \end{split}
\end{equation}
Let $f\in T_0^{m, p}(V)$ and $g\in T_1^{n, q}(V)$, we define
$\{f, g\}\in T_0^{m+n-1, p+q}(V)$ as 
\begin{equation*}
  \begin{split}
    \{f, g\}:=&\sum_{i=1}^{m-1}(-1)^{i(q-n-1)}(f\otimes
    \id^{\otimes q})(\id^{\otimes i-1}\otimes g\otimes \id^{m-i})-\\
    & \sum_{i=1}^{p-2}(-1)^{(q-n-1)(m-p+i)}(\id^{\otimes i}\otimes g\otimes
    \id^{\otimes p-i-1})(f\otimes \id^{\otimes n-1}).
  \end{split}
\end{equation*}
and $    \{g, f\}:=-(-1)^{(m-p-1)(n-q-1)} \{f, g\}.$
Define $\{f, g\}=0$ if $f\in T_0^{m, p}(V)$ and $g\in T_0^{n, q}(V)$.
Clearly, $\{\cdot,\cdot\}$ is skew-symmetry: 
For $f\in T^{m, p}(V)$ and $g\in T^{n, q}(V)$,
$$\{f, g\}=-(-1)^{(m-p-1)(n-q-1)}\{g, f\}.$$
By graphic computations we obtain the Jacobi identity, thus $(T^{*, *}(V), \{\cdot, \cdot\})$ is a graded Lie algebra. 

\begin{rem}
Note that $(T_1^{1, 0}(V), \{\cdot, \cdot\})$ is isomorphic to the Lie algebra $\mathfrak{gl}(V)$ with the usual Lie bracket, and $T_0^{m, p}(V)$ is the  representation of $\mathfrak{gl}(V)$ isomorphic to $V^{\otimes  p}\otimes (V^*)^{\otimes m}$. 
\end{rem}

For $m, p\in\Z_{\geq 0}$, we consider the embedding,
 \begin{equation*}
\begin{tabular}{ccccc}
$\theta_{m, p}$: &$T^{m,p}(V)$& $\hookrightarrow$ & $T^{m+1, p+1}(V)$\\
 &$(f, g)$&$\mapsto$&$(f\otimes \id_V, g\otimes \id_V)$\\
\end{tabular}
\end{equation*}
and the $k$-linear map
 \begin{equation*}
\begin{tabular}{ccccc}
$\phi_{m-1, p}$: &$T_0^{m-1,p}(V)$& $\rightarrow$ & $T_1^{m, p}(V)$\\
 &$f$&$\mapsto$&$\id_V\otimes f+(-1)^{p+m}f\otimes \id_V.$\\
\end{tabular}
\end{equation*}
Then for any $p\in \mathbb Z_{\geq 0}$, we have a complex $T^{*, p}(V)$:
\begin{equation}\label{equa-com}
\xymatrix@C=3pc{
 \cdots\ar[r] &  T^{m,p}(V) \ar[r]^-{\left( \begin{smallmatrix} 0 & 0\\ \phi_{m,p} &0 \end{smallmatrix}\right)}& T^{m+1, p}(V)\ar[r]^-{\left( \begin{smallmatrix} 0 & 0\\ \phi_{m+1,p} &0 \end{smallmatrix}\right)}   & T^{m+2, p}(V)\ar[r] &\cdots
  }
\end{equation}
and the embeddings $\theta$ induce embeddings of complexes 
$$T^{*, 0}(V) \xrightarrow{\theta_{*, 0}} T^{*, 1}(V)\xrightarrow{\theta_{*, 1}}\cdots \rightarrow T^{*, p}(V)\xrightarrow{\theta_{*, p}} \cdots.$$ 
Let us denote by $T_{\sg}^*(V)$ the colimit of the above inductive system. 
\begin{rem}
If $m\neq p$, then $\phi_{m, p}$ is injective and if $m=p$, $\Ker(\phi_{m,m})$
is a one-dimensional $k$-vector space with a basis $\{\id_{V^{\otimes m}}\}$.
As a consequence, the cohomology group  of the complex $T^{*, p}(V)$ is  \begin{equation}\label{equ-abc}
K^{m, p}(V):=
  \begin{cases}
   \frac{ T_1^{m, p}(V)}{\Imm(\phi_{m-1, p})} &
   \mbox{if $m\neq p$};\\
   k\id_{V^{\otimes m}}\oplus
   \frac{T_1^{m, m}(V)}{\Imm(\phi_{m-1,m})} & \mbox{if $m=p$}.
  \end{cases}
\end{equation}
\end{rem}
%\begin{lemma}\label{lemma5.4}
% For $m, p\in \Z_{\geq 0}$,  $\theta_{m, p}$ induces a morphism (still denoted by $\theta_{m,p}$),
%  $$
 %  \theta_{m, p}: K^{m, p}(V)\rightarrow K^{m+1, p+1}(V).$$
%\end{lemma}
%\pf Straightforward. \epf
\begin{rem}\label{rem-abc}
We have an inductive system
  \begin{equation*}
    \xymatrix@C=3pc{
    \cdots \ar[r] & K^{m, p}(V)\ar[r]^-{\theta_{m, p}} & K^{m+1, p+1}(V)\ar[r]^-{\theta_{
    m+1, p+1}} & \cdots
    }
  \end{equation*}
  Let us denote the colimit of this inductive system by
  $K_{\sg}^{m-p}(V)$ and we denote
  $$K_{\sg}^*(V):=\bigoplus_{n\in \Z} K_{\sg}^n(V).$$
  It is clear that $H^i(T_{\sg}^*(V)) \cong K_{\sg}^i(V)$ and $\theta_{m, m}(\id_{V^{\otimes m}})=\id_{V^{\otimes m+1}}$ for any $i\in \Z$ and $m\in \Z_{\geq 0}$. We denote the colimit of $\id_{V^{\otimes m}}$ by $\id_{V^{\otimes \infty}}$ in $K^0_{\sg}(V). $
\end{rem}

Now let us go back to the $r$-multiplied loops quiver $Q$. Recall that $A=kQ_0\oplus kQ_1$.
 % biggest distance of arrow to node. Yarou can use pt or cm as well.

Let $V:=k^r$ be the $k$-linear vector space of dimension $r$ with the canonical basis $e_1, \cdots, e_r\in k^r$. We can identify $k(Q_m//Q_p)$ with $\Hom_k(V^{\otimes m}, V^{\otimes p})$
as follows,
\begin{equation}\label{formula-cor}
\begin{tabular}{ccccc}
 &$k(Q_m//Q_p)$& $\rightarrow$ & $\Hom_k(V^{\otimes m}, V^{\otimes p})$\\
 &$(x_1\cdots x_m, y_1\cdots y_p)$&$\mapsto$&$\delta_{e_{x_1}\otimes \cdots\otimes e_{x_m}, e_{y_1}\otimes \cdots \otimes e_{y_p}}$\\
\end{tabular}
\end{equation}
where $x_i, y_j\in\{1, \cdots, r\}$ and $\delta_{e_{x_1}\otimes \cdots\otimes e_{x_m}, e_{y_1}\otimes \cdots \otimes e_{y_p}}$ is defined as follows,
\begin{equation*}
  \delta_{e_{x_1}\otimes \cdots\otimes e_{x_m}, e_{y_1}\otimes \cdots \otimes e_{y_p}}(e_{x_1'}\otimes \cdots\otimes e_{x_m'})=
  \begin{cases}
   e_{y_1}\otimes \cdots \otimes e_{y_p}  & \mbox{if $x_i=x_i'$ for all $1\leq i \leq m$};\\
    0 &\mbox{othewise}.
  \end{cases}
\end{equation*}
Thus we have the following isomorphism for $m, p\in \Z_{\geq 0}$,
\begin{equation*}
  \begin{tabular}{ccccc}
 $F_{m,p}:$ &$k(Q_m//Q_p)\oplus k(Q_m//Q_{p+1})$& $\xrightarrow{\cong}$ & $T^{m, p}(V).$\\
 \end{tabular}
\end{equation*}
%Then we have the following proposition.
\begin{thm}\label{prop-prop-two}
Let $Q$ be the $r$-multiplied loops quiver and $A$ be its radical square zero algebra over a field $k$ of characteristic zero. Then
\begin{enumerate}
  \item for any  $p\in\Z_{\geq0}$, we have an isomorphism of complexes
  $$F_{*, p}: C^*(r, \Omega^p_{nc}(A))\xrightarrow{\cong} T^{*, p}(V).$$
%  \begin{equation*}
%    \xymatrix@C=2.5pc{
%   \cdots\ar[r] & k(Q_m//Q_p)\oplus k(Q_m//Q_{p+1})\ar[d]^{F_{m, p}} \ar[r]^-{\left( \begin{smallmatrix} 0 & 0\\ D_{m,p} &0 \end{smallmatrix}\right)}& k(Q_{m+1}//Q_p)\oplus k(Q_{m+1}//Q_{p+1})\ar[d]^{F_{m+1, p}}\ar[r] &\cdots\\
%   \cdots\ar[r] & T^{m,p}(V) \ar[r]^-{\left( \begin{smallmatrix} 0 & 0\\ \phi_{m,p} &0 \end{smallmatrix}\right)}& T^{m+1, p}(V)\ar[r] & \cdots
 %   }
 % \end{equation*}
%  where we recall that $D_{m,p}$ is defined in Proposition \ref{prop-new}. 
This induces an isomorphism  $F_{*, p}: \HH^{*}(A, \Omega_{nc}^p(A))\rightarrow K^{*, p}(V)$.
  \item The following diagram commutes, for $m, p\in \Z_{>0}$.
  \begin{equation*}
    \xymatrix{
    \HH^{m}(A, \Omega_{nc}^p(A)) \ar[d]^-{\cong}_-{F_{m, p}}\ar[r]^-{\theta_{m, p}} & \HH^{m+1}(A, \Omega_{nc}^{p+1}(A))\ar[d]_{\cong}^-{F_{m+1, p+1}}\\
    K^{m, p}(V)\ar[r]^-{\phi_{m, p}}& K^{m+1, p+1}(V)
    }
  \end{equation*}
  %As a consequence, $\phi_{m, p}$ is injective since $\
  \item The following diagram commutes for $m, n\in \Z_{>0}$ and  $ p, q \in \Z_{\geq 0}$
  \begin{equation*}
    \xymatrix{
     \Hom_{E-E}(r^{\otimes_E m}, \Omega_{nc}^p(A))\times \Hom_{E-E}(r^{\otimes_E n}, \Omega_{nc}^q(A)) \ar[d]^-{F_{m, p}\times F_{n, q}}_-{\cong}\ar[r]^-{[\cdot,\cdot]} &\Hom_{E-E}(r^{\otimes_E m+n-1}, \Omega_{nc}^{p+q}(A))\ar[d]^-{F_{m+n-1, p+q}}_-{\cong}\\
    T^{m, p}(V) \times T^{n, q}(V) \ar[r]^{\{\cdot,\cdot\}} &  T^{m+n-1, p+q}(V)
    }
  \end{equation*}
  where the Gerstenhaber bracket $[\cdot,\cdot]$ is defined in (\ref{defn-bracket}) and $\{\cdot, \cdot\}$
  is defined in (\ref{defn-bracket-T}). Here we identify $$\Hom_{E-E}(r^{\otimes_Em}, \Omega_{nc}^p(A))$$ with
  $$k(Q_m//Q_p)\oplus k(Q_m//Q_{p+1})$$ by
  Lemma \ref{lemma-new}.
\end{enumerate}
\end{thm}
\pf  Assertions (1) and (2) follow from straightforward computations. Let us
verify the Assertion (3). Let $$(x_1 \cdots x_m, y_1 \cdots y_{p+1})
\in \Hom_{E-E}(r^{\otimes_E m}, \Omega_{nc}^p(A))$$ and $$(a_1 \cdots
a_n, b_1\cdots b_{q+1})\in \Hom_{E-E}(r^{\otimes_E n}, \Omega_{nc}^q(A)).$$
By (\ref{formula-cor}), these two elements correspond to 
$$\delta_{e_{x_1}\otimes \cdots \otimes e_{x_m}, e_{y_1}\otimes
\cdots \otimes e_{y_{p+1}}}\in T^{m ,p }$$ and $$ \delta_{e_{a_1}\otimes \cdots \otimes e_{a_{n}}, e_{b_1}\otimes
\cdots \otimes e_{b_{q+1}}}\in T^{n, q},$$ respectively.  Then under this correspondence, we observe that Formulae (\ref{defn-bracket}) and (\ref{defn-bracket-T}) coincide. This proves the theorem. \epf

%Therefore, as a consequence, we obtain the following corollary.
%\begin{cor}\label{cor-last}
 % Let $V$ be a finite dimensional vector space of dimensional $r$ over a field $k$ of characteristic zero. Then
 % \begin{enumerate}
 % \item $T^{*,*}(V)[1]$ is a differential graded Lie algebra, where the grading
%  is defined as follows, for $n\in \Z$,
%  $$T^{*, *}(V)_n:=\bigoplus_{\substack{m, p\in\Z_{\geq 0}\\ m-p=n}} T^{m , p}(V)$$
 % and the differential is induced by
  %$\xymatrix@C=3pc{
 % T^{m,p}(V) \ar[r]^-{\left( \begin{smallmatrix} 0 & 0\\ \phi_{m,p} &0 \end{smallmatrix}\right)}& T^{m+1, p}(V),
%  }$
%  In particular, $K^{*, *}(V)$ is a $\Z$-graded Lie algebra with the induced
%  grading and induced Lie bracket;
%  \item The embeddings  $\theta_{m, p}:T^{m, p}(V)\rightarrow T^{m+1, p+1}(V)$ respect to the Lie bracket, namely,  $\theta([f, g])=[\theta(f), g]$. As a consequence,  $K^*_{\sg}(V)$ is  isomorphic to $\HH_{\sg}^*(A, A)$ as graded Lie algebras.
%  \end{enumerate}
%\end{cor}
%\pf This is an immediate consequence of Proposition \ref{prop-prop-two}.
%\epf

Recall that in Remark \ref{rem-abc} we denote by $\id_{V^{\otimes \infty}}$, the colimit of $\id_{V^{\otimes m}}$ in $K_{\sg}^0(V)$. Note that from (\ref{equ-abc}) we have a canonical projection 
$$\epsilon: K_{\sg}^*(V) \twoheadrightarrow k$$ sending $\id_{V^{\otimes \infty}}$ to $1$. 
Then there is a unique $k$-algebra $(K_{\sg}^*(V), \cdot)$ such that $\epsilon$ is an algebra morphism and that $\Ker(\epsilon)$ has trivial multiplication. In other words, $K_{\sg}^*(V)$ is the trivial extension of the algebra $k$ and the $k$-module $\Ker(\epsilon)$. 

\begin{cor}\label{cor-last}
Let $Q$ be the $r$-multiplied loops quiver and $A$ be its radical square zero algebra over a field $k$ of characteristic zero. Then $(\HH_{\sg}^*(A, A), \cup, [\cdot, \cdot])$ is isomorphic to $(K_{\sg}^*(V), \cdot, \{\cdot, \cdot\})$ as Gerstenhaber algebras. 
\end{cor}
\pf
From Remark  \ref{rem-abc} and Theorem \ref{prop-prop-two},  it follows that $\id_{V^{\otimes \infty}}$ represents the unit $1\in \HH^0_{\sg}(A, A)$. %From (\ref{equ-abc}) and Theorem \ref{prop-prop-two},  we have an augmentation $$\epsilon: \HH_{\sg}^*(A, A)\twoheadrightarrow k\id_{V^{\otimes \infty}}.$$
Thus from Corollary \ref{cor-cup}  and its proof, we infer that the ideal $\Ker(\epsilon)$ has trivial multiplication. Therefore this corollary follows from Theorem \ref{prop-prop-two}. \epf
\begin{rem}\label{rem-conter1}
Since the cup product $\cup$ of $\HH_{\sg}^*(A, A)$ vanishes except in degree zero and  the Lie bracket $[\cdot, \cdot]$ does not.  Thus it is impossible to endow the Gerstenhaber algebras $(\HH_{\sg}^*(A, A), \cup, [\cdot,\cdot])$ and $(\HH^*(A, A), \cup, [\cdot, \cdot])$ with BV algebra structures.
\end{rem}

%At the end of this section, let us describe the Lie algebra structure on $\HH_{\sg}^1(A, A)$, where $A$ is the radical square zero algebra of the two loops quiver.  From Corollary \ref{cor-last},

\bibliographystyle{plain}

\end{document}